\documentclass[11pt]{article}
\usepackage[margin=1in]{geometry}
\usepackage[utf8]{inputenc}
\usepackage[english]{babel}
\usepackage{amsmath,amsthm}
\usepackage{color}
\usepackage{amsfonts,enumitem}
\usepackage{amssymb}
\usepackage{hyperref}
\usepackage{cleveref}
\Crefname{algocf}{Algorithm}{Algorithms}
\usepackage{graphicx,tikz}
\usepackage{dsfont}

\usepackage{algorithm}
\usepackage{algorithmic}

\newcommand{\crit}{\mathrm{crit}}
\newcommand{\conv}{\mathrm{conv}}
\newcommand{\dist}{\mathrm{dist}}
\newcommand{\closedball}{\overline{\mathbb{B}}}

\newcommand{\partialc}{\partial^c}
\newcommand{\supp}{\operatorname{supp}}

\newcommand*{\SET}[1]  {\ensuremath{\mathbb{#1}}}
\newcommand{\argmax}{\operatornamewithlimits{argmax}}
\newcommand{\argmin}{\operatornamewithlimits{argmin}}
\newcommand{\R}{\SET{R}}
\newcommand{\E}{\SET{E}}
\newcommand{\D}{\SET{D}}
\newcommand{\N}{\SET{N}}
\newcommand{\KL}{\operatorname{KL}}

\newcommand{\Graph}{\operatorname{Graph}}

\newcommand{\di}{\text{d}}

\newcommand{\diff}{\operatorname{diff}}

\newtheorem{theorem}{Theorem}[section]
\newtheorem{corollary}{Corollary}[section]
\newtheorem{lemma}{Lemma}[section]
\newtheorem{proposition}{Proposition}[section]
\theoremstyle{definition}
\newtheorem{definition}{Definition}[section]
\newtheorem{assumption}{Assumption}

\begin{document}
\title{Entropic smoothing and gradient approximation for Wasserstein distributionally robust optimization}
\author{Tam Le \\
        Universit\'{e} Paris Cit\'{e} \\
        Paris, 75013 FRANCE \\
        \href{mailto:tamle@lpsm.paris}{\texttt{tamle@lpsm.paris}}}
\date{}
\maketitle

\begin{abstract}
Wasserstein distributionally robust optimization offers a framework for model fitting in machine learning under potential shifts in the data distribution. We study a regularized variant of this problem in which entropic smoothing produces a sampled approximation of the original objective. We establish convergence of the approximate gradients to subgradients of the unregularized objective as the regularization parameter vanishes, enabling convergence guarantees for stochastic gradient methods.
We obtain qualitative convergence results under general assumptions and provide convergence rates under additional regularity. In particular, we prove rates for the convergence of the unregularized objective values, up to sampling errors, when the regularization level is decreased across iterations. Our analysis also yields results of independent interest on approximating
the subdifferentials of maximum functions by gradients of their smooth
approximations, together with growth-based lower bounds on optimal dual
variables in Wasserstein distributionally robust optimization.
\end{abstract}



\section{Introduction}

Machine learning models are commonly trained by empirical risk minimization. Given a loss function $\ell$, this amounts to solving
\begin{equation*}
    \min_{\theta \in \Theta} 
    \frac{1}{n} \sum_{i=1}^n \ell(\theta,\xi_i)
    =
    \min_{\theta \in \Theta}
    \E_{\xi \sim \widehat{P}_n}[\ell(\theta,\xi)],
\end{equation*}
where $\Theta \subset \R^p$ is the parameter space and $\widehat{P}_n$ is the empirical distribution associated with the data sample $\xi_1,\ldots,\xi_n$ taking values in a space $\mathcal{Z}$. Distributionally robust optimization has been widely used in this context to improve the robustness of learned models with respect to distributional uncertainty \cite{blanchet2025distributionally}. In this framework, one minimizes the worst-case risk over an ambiguity set of probability distributions close to $\widehat{P}_n$. Among the possible choices, ambiguity sets defined by optimal transport costs \cite{peyre2019comput,villani} have emerged as a prominent class \cite{logisticregWDRO2016,MohajerinEsfahaniKuhn2018,kuhn2019wasserstein}, supported by attractive statistical properties \cite{gao2022variationregularization,azizian2023exact,le2025universal,gaofinitesample2022}.

Let $c:\mathcal{Z}\times \mathcal{Z}\to \R_+$ be a cost function. For two probability distributions $P$ and $Q$ on $\mathcal{Z}$, the corresponding optimal transport cost is defined by
\begin{equation*}
    W_c(P,Q)
    =
    \inf_{\substack{\pi \in \mathcal{P}(\mathcal{Z}\times \mathcal{Z})\\
    \pi_1=P,\ \pi_2=Q}}
    \E_{(\xi,\zeta)\sim \pi}[c(\xi,\zeta)],
\end{equation*}
where $\mathcal{P}(\mathcal{Z}\times \mathcal{Z})$ denotes the set of probability measures on $\mathcal{Z}\times \mathcal{Z}$, and $\pi_1,\pi_2$ denote the first and second marginals of $\pi$. When $c$ is the $q$-th power of a distance, $W_c^{1/q}$ recovers the usual $q$-Wasserstein distance. Wasserstein distributionally robust optimization (WDRO) then considers problems of the form
\begin{equation}
    \label{eq:wdro}
    \min_{\theta \in \Theta}
    \max_{\substack{Q \in \mathcal{P}(\mathcal{Z})\\
    W_c(\widehat{P}_n,Q) \leq \rho}}
    \E_{\xi \sim Q}[\ell(\theta,\xi)],
\end{equation}
where $\rho>0$ represents a fixed uncertainty level. A key feature of Wasserstein ambiguity sets is that the inner worst-case expectation admits a dual representation \cite{blanchet2019quantifying,gao2023wassersteindistance,zhang2024short} as
\begin{equation}
    \label{eq:not_joint_objective}
    \max_{\substack{Q \in \mathcal{P}(\mathcal{Z})\\
    W_c(\widehat{P}_n,Q) \leq \rho}}
    \E_{\xi \sim Q}[\ell(\theta,\xi)]
    =
    \min_{\lambda \geq 0}
    \E_{\xi \sim \widehat{P}_n}
    \left[
        \sup_{z \in \mathcal{Z}}
        \big\{
            \ell(\theta,z) + \lambda(\rho-c(\xi,z))
        \big\}
    \right].
\end{equation}
Consequently, the min--max problem \eqref{eq:wdro} can be reformulated as the joint minimization problem
\begin{equation}
    \label{eq:joint_objective_unreg}
    \min_{\theta \in \Theta,\ \lambda \geq 0}
    F(\theta,\lambda)
    :=
    \E_{\xi \sim \widehat{P}_n}
    \left[
        \sup_{z \in \mathcal{Z}}
        \big\{
            \ell(\theta,z) + \lambda(\rho-c(\xi,z))
        \big\}
    \right].
\end{equation}
Despite this reformulation, the supremum in \eqref{eq:joint_objective_unreg} is generally difficult to compute, especially when the function inside the supremum is nonconcave in $z$. To address this difficulty, regularization and smoothing techniques have been proposed \cite{wang2023sinkhorn,azizian2023regularization}. A common approach replaces the supremum by a smooth entropic approximation. Given a reference distribution $\pi_0$ on $\mathcal{Z}$ and a regularization parameter $\beta>0$, one considers, for each data point $\xi$, the smoothed quantity
\begin{equation}
    \label{eq:entropic_smoothing_phi_beta}
    \beta \log
    \E_{z \sim \pi_0}
    \left[
        \exp\left(
            \frac{
                \ell(\theta,z) + \lambda(\rho-c(\xi,z))
            }{\beta}
        \right)
    \right].
\end{equation}
The reference distribution $\pi_0$ is chosen so that samples can be easily generated. This leads to the smoothed objective
\begin{equation}
    \label{eq:regularized_objective}
    F^\beta(\theta,\lambda)
    :=
    \frac{1}{n}\sum_{i=1}^n
    \beta \log
    \E_{z \sim \pi_0}
    \left[
        \exp\left(
            \frac{
                \ell(\theta,z) + \lambda(\rho-c(\xi_i,z))
            }{\beta}
        \right)
    \right].
\end{equation}
In practice, the expectation with respect to $\pi_0$ can be approximated by sampling. More precisely, given an i.i.d. sample $z_1,\ldots,z_m \sim \pi_0$, this yields the finite-sample approximation
\begin{equation}
    \label{eq:regularized_objective_samples}
    F^{\beta,m}(\theta,\lambda)
    :=
    \frac{1}{n}\sum_{i=1}^n
    \beta \log
    \left(
        \frac{1}{m}
        \sum_{j=1}^m
        \exp\left(
            \frac{
                \ell(\theta,z_j) + \lambda(\rho-c(\xi_i,z_j))
            }{\beta}
        \right)
    \right).
\end{equation}
Standard stochastic-gradient methods can then be applied to minimize \eqref{eq:regularized_objective_samples}. This principle underlies recent implementations of WDRO \cite{wang2023sinkhorn,vincent2024skwdro,liu2025dro}.

\paragraph{Contributions.}
In this work, we study the approximation properties of the regularized and sampled objectives $F^\beta$ and $F^{\beta,m}$ as the regularization parameter $\beta$ vanishes and the number of approximation samples $m$ grows. In \Cref{subsection:qualitative_convergence}, we establish convergence results for the gradients of the regularized objectives toward the Clarke subdifferential \cite{clarke1990optimization} of the nonsmooth unregularized objective $F$ \eqref{eq:joint_objective_unreg}. More precisely, uniform convergence holds on compact sets for which the dual variable is bounded away from zero. Specifically, we show that for every $\varepsilon>0$, if $\beta>0$ is sufficiently small and $m$ is sufficiently large, then
\begin{equation}
\label{eq:contribution_convergence}
\nabla F^{\beta,m}(\theta,\lambda)
\in
\partialc F\big((\theta,\lambda)+\varepsilon\bar{\mathbb{B}}\big)
+ \varepsilon\bar{\mathbb{B}}
\quad
\text{for all }(\theta,\lambda)\in \Theta\times\Lambda,
\end{equation}
where $\Lambda \subset (0,\infty)$. We propose conditions under which one may restrict attention to such sets $\Lambda$ in \Cref{subsection:truncation}. Under additional regularity assumptions, we derive explicit approximation estimates in \Cref{subsection:quantitative_estimates}. 
 
 We use these approximation results to analyze stochastic-gradient methods applied to the smoothed sampled objectives. In \Cref{subsection:gradient_methods}, we present algorithms with vanishing regularization parameters $\beta_k \to 0$ across iterations and their convergence guarantees for minimizing the nonsmooth objective.

\subsection{Connections with related literature}
\paragraph{Entropic smoothing.}

Regularization techniques for WDRO were introduced in
\cite{wang2023sinkhorn,azizian2023regularization}. In particular, the cost
regularization proposed in \cite{azizian2023regularization} corresponds to
the entropic smoothing in \eqref{eq:entropic_smoothing_phi_beta}. Algorithms for convex objectives were also proposed in \cite{wang2023sinkhorn}. Optimization guarantees for the joint objectives
\eqref{eq:joint_objective_unreg} and 
\eqref{eq:regularized_objective_samples} were studied in several works \cite{liu2025stochastic,piat2026distributionally}. More broadly, entropic smoothing is a well-established technique in variational analysis that provides a natural smooth approximation to maximum-value functions \cite{fang1996solving,li1997entropic,lin2014solving,alcantara2024theoretical}.  We
highlight below the main contributions of our work relative to this
literature on entropic smoothing and distributionally robust optimization.

\begin{itemize}

	\item We provide a self-contained proof of the graphical convergence of the smoothed gradient to the subdifferential under general assumptions, notably allowing for unbounded sample spaces. This convergence ensures that, for sufficiently small fixed regularization, critical points of the smoothed problem remain close to those of the original nonsmooth problem. A uniform gradient-approximation bound is stated in   \cite{piat2026distributionally}. Such a bound cannot hold in general for a nonsmooth limiting objective, due to transitions between active maximizers. The proof in \cite{piat2026distributionally} actually establishes pointwise convergence, whereas our graphical convergence provides the appropriate uniform notion in this setting. Pointwise gradient convergence is also established in
\cite{liu2025stochastic} and used to derive a convergence result when considering a sequence of subproblems with decreasing smoothing parameters. Our uniform graphical convergence ensures the consistency of critical points for a fixed, sufficiently small regularization parameter, corresponding to a more practical setting.
	
\item We treat unbounded sample spaces. This is done by adapting standard
uniform-integrability and tightness arguments for Gibbs measures. We believe their use in this setting may be of independent interest for variational problems involving similar sampling approximations. By comparison, previous studies consider compact sample spaces \cite{liu2025stochastic,piat2026distributionally} or sampling distributions with growing compact supports \cite{alcantara2024theoretical}.

\item We derive optimization bounds for regularization schedules that may
vary across iterations, including vanishing regularization schedules. This provides a basis for future algorithmic developments in which the regularization level and the sampling accuracy  would be adapted jointly. This is inspired by recent developments in entropic optimal transport \cite{chizat2024annealed,genans2025decreasing}.
\end{itemize}

\paragraph{Function to subdifferential convergence.} Another possible route to subdifferential convergence would be to use
rigidity results under which convergence of objective functions implies
convergence of their subdifferentials; see
\cite{attouch2006convergence,levy1995partial,schechtman2024gradient}. Under
our assumptions, we expect that such a connection can be established
through uniform semidifferentiability arguments of the type developed in
\cite{zolezzi1985continuity}. This connection is nevertheless not direct
and would require further development. Our analysis makes explicit
the convergence mechanism at work by combining convergence of the Gibbs
measures with standard notions from variational analysis, including the
 envelope formula \cite[Th. 2.8.2]{clarke1990optimization} and graphical convergence. Under weak convexity,
the connection is more direct: approximation bounds for objective functions
can be converted into approximation bounds for their subdifferentials; see
\cite{davis2022graphical}. We use this principle to derive the quantitative
estimates in \Cref{th:quantitative_approx_graph}.

\paragraph{Dual solution lower bound.} Applying our results to the unrestricted WDRO problem requires a lower
bound $\lambda_{\min}>0$ such that the restriction
$\lambda\geq\lambda_{\min}$ preserves the optimal value and satisfies
our growth assumptions. Related lower bound conditions have been used
to study statistical properties of robust objectives
\cite{azizian2023exact,le2025universal,gaofinitesample2022}.
In \Cref{subsection:truncation}, we discuss the relationship between
lower bounds on optimal dual variables and the growth of the loss
relative to the transport cost.

\paragraph{Sampling and numerical study.}
Our work is primarily theoretical, focusing on the approximation effects of entropic regularization and finite sampling. In particular, we establish general qualitative consistency results for finite-sample approximations based on i.i.d.\ samples, without addressing the design or computational performance of advanced sampling schemes. For practical sampling strategies and numerical implementations, we refer to \cite{wang2023sinkhorn,vincent2024skwdro,liu2025dro}, and for numerical experiments and quantitative sampling bounds on compact sets, to \cite{piat2026distributionally,liu2025stochastic}.

\subsection{Organization of the paper}
The remainder of the paper is organized as follows. \Cref{sec:main_results} presents the main results, including the qualitative approximation theorem, quantitative estimates, and convergence guarantees for gradient-based methods. The subsequent sections provide the proofs and auxiliary results needed for the analysis. Specifically, \Cref{section:convergence_entropic_general} establishes convergence results for gradients of entropic regularization, while \Cref{sec:convergence_rates_general} develops abstract convergence rates for approximate stochastic gradient methods. Finally, \Cref{sec:proofs} contains the proofs of the main results, which connect the WDRO setting to our general results.
\bigskip

\paragraph{Notations and essential definitions} For $d \geq 1$, let  $I_d$ denote the identity matrix of size $d$.  $\|\cdot\|$ denotes the Euclidean norm and $\closedball$ is the closed unit ball. $\closedball(x,r)$ is the closed ball centered at $x$ and radius $r\geq 0$. For $q \geq 1$, when $A \subset \R^q$, we denote the closure of $A$ by $\Bar{A}$ and its complement by $A^c$.  We denote the distance of $x \in \R^q$ from a subset $A$ by  $\dist(x, A) := \inf_{a \in A} \|x - a\|$, and when $A, B$ are nonempty subsets, we define the excess distance as 
\begin{equation}
    \label{eq:excess_distance}
    \D \left(A,B \right) = \sup_{a \in A} \dist(a, B).
\end{equation}
and the Hausdorff distance is $\mathbb{H}(A,B) = \max(\D(A,B), \D(B,A))$, allowing the value $+\infty$.

-- \emph{Subdifferentials and set-valued maps.}
Let $f:\mathbb{R}^q\to\mathbb{R}$ be locally Lipschitz, and denote by
$\diff f$ the set of points at which $f$ is differentiable. In particular, $\diff f$ has full Lebesgue measure. The
\emph{Clarke subdifferential} of $f$ at $w\in\mathbb{R}^q$ is defined by
\begin{equation}
	\label{eq:clarke}
	\partialc f(w)
	:=
	\conv\left\{
	v\in\mathbb{R}^q
	\ \middle|\
	\begin{array}{l}
		\text{there exists }(w_k)_{k\in\mathbb{N}}\subset\diff f
		\text{ such that}\\
		w_k\to w
		\text{ and }
		\nabla f(w_k)\to v
	\end{array}
	\right\},
\end{equation}
where $\conv$ denotes the convex hull. When $K$ is compact and convex, we define the \emph{critical points} restricted to $K$ as $\crit_{K} f := \{w \in K \ : \ 0 \in \partialc f(w) + N_K(w)\}$
where $N_K(w)$ is the normal cone to $K$ at $w$. $\Pi_{K}$ is the projection onto $K$. The graph of a set-valued map $D : K \subset \R^p \rightrightarrows \R^p$ is defined as 
\begin{equation}
	\label{eq:graph}
	\Graph_K D := \left\{(x,y) \ : \ x \in K, y \in D(x)\right\}.
\end{equation}

-- \emph{Level sets.} For a function $h: \R^q \to \R$ and $a \in \R$, we denote sublevel sets by $\{h \leq a\} := \left\{x \in \R^q \ : \ h(x) \leq a\right\}$ and we use analogous shorthand for superlevel sets and their counterparts with strict inequalities. If $\mathcal{Z}$ is a metric space, we call a continuous function $\psi : \mathcal{Z} \to \R$ \textit{coercive} when $\{\psi \leq a\}$ is compact for all $a \in \R$.

\emph{Probability spaces.} Throughout the paper $\mathcal{Z}$ is a Polish space\footnote{complete metric space with a countable dense subset}.
We write $\mathbb{P}(A)$ and $\mathbf{1}_A$ for the probability and indicator of an event $A$, respectively. For a measurable space $\mathcal{Z}$, $\mathcal{P}(\mathcal{Z})$ denotes the set of probability measures on $\mathcal{Z}$. For $z \in \mathcal{Z}$, $\delta_z$ denotes the Dirac measure at $z$, and for $\pi \in \mathcal{P}(\mathcal{Z})$, $\operatorname{supp}(\pi)$ denotes its support.

\section{Main results}
\label{sec:main_results}

In this part, we present the main results on the regularized WDRO problem introduced above.  Throughout this part we work on a compact parameter set and restrict the dual variable to a compact interval bounded away from zero. We discuss when minimizing $F$ over $\Theta \times \Lambda$ preserves the optimal value in \Cref{subsection:truncation}, see also \cite{le2025universal,azizian2023exact} for compact sample spaces. Thus let $\Theta \subset \R^p$ be compact and convex and let $\Lambda = [\lambda_{\min}, \lambda_{\max}]$, where
$0 < \lambda_{\min} < \lambda_{\max}<\infty$. 

\subsection{Qualitative convergence}
\label{subsection:qualitative_convergence}

We consider the WDRO objective  for $(\theta,\lambda)\in \Theta\times\Lambda$,
\begin{equation}
\label{eq:WDRO_objective_main}
    F(\theta,\lambda)
:=
\frac1n\sum_{i=1}^n \sup_{z\in\mathcal Z}
\Big\{\ell(\theta,z)+\lambda\bigl(\rho-c(\xi_i,z)\bigr)\Big\}.
\end{equation} Given a  regularization parameter $\beta>0$, the corresponding entropic smoothing is given by
\begin{equation}
\label{eq:regularized_WDRO_main}
    F^{\beta}(\theta,\lambda)
    :=
    \frac{1}{n} \sum_{i=1}^n
    H_i^\beta(\theta,\lambda),
\end{equation}
where $$H_i^\beta(\theta,\lambda):= \beta \log \E_{z \sim \pi_0}
\left[
\exp \left(
\frac{\ell(\theta,z)+\lambda(\rho-c(\xi_i,z))}{\beta}
\right)
\right]$$   
and $\pi_0$ is a reference distribution on $\mathcal{Z}$.

Our first result shows that, as the regularization  level $\beta$ vanishes, the gradients of the smooth objectives $F^\beta$ approximate the Clarke subdifferential \eqref{eq:clarke} of the nonsmooth objective $F$ uniformly on $\Theta\times\Lambda$. In particular, this implies convergence of the restricted critical sets. We work under the following set of assumptions. 

\begin{assumption}
\label{ass:regularity_ell_c_merged}
\begin{enumerate}
    \item[]
    \item \(\supp(\pi_0)=\mathcal Z\).

    \item \(\ell\) and $c$ are continuous. For every \(z\in\mathcal Z\), the map \(\theta\mapsto \ell(\theta,z)\) is differentiable and \(\nabla_\theta \ell\) is jointly continuous on \(\Theta\times\mathcal Z\).

    \item For all $\xi \in \mathcal{Z}$, $c(\xi, \cdot)$ is non-negative, coercive on $\mathcal{Z}$ and   $\pi_0$-integrable. Furthermore, there exists  a $\pi_0$-integrable function \(\varphi:\mathcal Z\to[0,\infty)\), and constants $L \geq 0$, \(\bar{\lambda}\in[0,\lambda_{\min})\) such that for all  $i=1,\dots,n$ and \((\theta,z)\in\Theta\times\mathcal{Z}\),
    \[
    |\ell(\theta,z)|\le L + \bar{\lambda} c(\xi_i, z) \quad\text{and}
    \quad
    \|\nabla_\theta \ell(\theta,z)\|\le \varphi(z).
    \]
\end{enumerate}
\end{assumption}

\begin{theorem}[Uniform gradient convergence] 
    \label{th:gradient_conv_ell_c}
Let $\varepsilon > 0$. Under \Cref{ass:regularity_ell_c_merged} there exists $\bar{\beta} > 0$ such that  for any $\beta \in (0, \bar{\beta}]$ and all  $(\theta, \lambda) \in \Theta \times \Lambda$,
\begin{enumerate}
    \item $\nabla F^{\beta}(\theta,\lambda) \in \partialc F((\theta,\lambda) + \varepsilon\bar{\mathbb{B}}) + \varepsilon \bar{\mathbb{B}}.$
    \item $\crit_{\Theta \times \Lambda} F^{\beta} \subset (\crit_{\Theta \times \Lambda} F) + \varepsilon \bar{\mathbb{B}}.$
\end{enumerate}
\end{theorem}
\begin{proof} See \Cref{subsection:convergence_gradient_proof}.
\end{proof}

\paragraph{On sampling approximation.}
Given $m$  i.i.d  draws from $\pi_0$, $z_1, \ldots, z_m$, we consider sampling approximations for $H_i^{\beta}$ and $F^{\beta}$ as
$$H_i^{\beta,m}(\theta,\lambda) = \beta \log \left( \frac{1}{m}\sum_{j=1}^m \exp  \left(
\frac{\ell(\theta,z_j)+\lambda(\rho-c(\xi_i,z_j))}{\beta}
\right) \right)$$ and $F^{\beta,m}(\theta,\lambda) = \frac{1}{n} \sum_{i=1}^n H_i^{\beta,m}(\theta,\lambda).$ Under the continuity and integrability conditions from \Cref{ass:regularity_ell_c_merged},  we have the following approximation guarantee.

\begin{lemma} \label{lem:LLN_H_beta_m} Under \Cref{ass:regularity_ell_c_merged}, for any $\beta,\varepsilon>0$, almost surely  there exists $m_{\beta,\varepsilon} \geq 1$ such that for any $m \geq m_{\beta,\varepsilon}$ and $i=1,\ldots,n$,
\begin{equation*}
       \sup_{\Theta \times \Lambda}|H_i^{\beta,m} - H_{i}^\beta|
       + \sup_{\Theta \times \Lambda} \|\nabla H_i^{\beta,m} - \nabla H_{i}^\beta\| \leq \varepsilon.
\end{equation*}
If, in addition, there exists a nonnegative $\pi_0$-integrable function $a$ such that
$\sup_{\theta\in\Theta}|\ell(\theta,z)|\leq a(z)$, then
\[
\mathbb E\!\left[\max_{1\leq i\leq n}\sup_{w\in K}
|H_i^{\beta,m}(w)-H_i^\beta(w)|\right]\longrightarrow0,
\qquad K=\Theta\times\Lambda.
\]
Consequently, for every $\beta,\varepsilon>0$, there exists a deterministic
$m_{\beta,\varepsilon}$ such that this expectation is at most
$\varepsilon$ for all $m\geq m_{\beta,\varepsilon}$.
\end{lemma}

\begin{proof} See \Cref{subsection:convergence_gradient_proof}.
\end{proof}
A direct consequence of \Cref{th:gradient_conv_ell_c} and  \Cref{lem:LLN_H_beta_m} is as follows.

\begin{theorem}[Uniform gradient convergence with sampling]
\label{th:gradient_consist_subsampling}
Let $\varepsilon > 0$. Under \Cref{ass:regularity_ell_c_merged}, there exists $\bar{\beta} > 0$ such that for any $\beta \in (0, \bar{\beta}]$, almost surely, there exists $\overline{m} \in \mathbb{N}$ such that for any $m \geq \overline{m}$ and all $(\theta, \lambda) \in \Theta \times \Lambda$,
\begin{enumerate}
    \item $\nabla F^{\beta,m}(\theta,\lambda) \in \partialc F((\theta,\lambda) + \varepsilon\bar{\mathbb{B}}) + \varepsilon \bar{\mathbb{B}}.$
    \item $\crit_{\Theta \times \Lambda} F^{\beta,m} \subset (\crit_{\Theta \times \Lambda} F) + \varepsilon \bar{\mathbb{B}}.$
\end{enumerate}
\end{theorem}

\subsection{Quantitative estimates for entropic smoothing}
\label{subsection:quantitative_estimates}

We now turn to quantitative estimates. In this subsection, we impose stronger regularity assumptions that give us more explicit control over the approximation error, in terms of the smoothing parameter $\beta > 0$.

We focus on the case of a quadratic transport cost and assume enough smoothness in the transported variable to obtain a uniform rate for the regularized objective. Under these regularity assumptions, this rate is then converted into a graphical estimate for the corresponding subdifferentials.
\begin{assumption}
	\label{ass:WDRO_smooth_quadratic}
	$\mathcal Z=\mathbb R^d$ and the following conditions hold.
	\begin{enumerate}
		
		\item The loss $\ell$ is jointly continuous on
		$\R^p\times\mathcal Z$ and $c$ is the quadratic cost, $c(z,z')=\|z-z'\|^2$ for all $z,z' \in \mathcal{Z}$. Moreover, for all $z \in \mathcal{Z}$, $\ell(\cdot,z)$ is differentiable and 
		$\nabla_\theta\ell$ is jointly continuous.
		
		\item There exists $\tau_\ell>0$ such that for all $z \in \mathcal{Z}$, $\ell(\cdot, z)$ is $\tau_\ell$-weakly convex.

		\item $\nabla_z \ell$ is jointly locally Lipschitz continuous on $\R^p\times\mathcal Z$.

		\item There exist a  continuous and $\pi_0$-integrable function
		\(\Psi:\mathcal Z\to[0,\infty)\), a locally bounded function
		\(C : \R^p \to [0,\infty) \), \(L\geq0\), and
		\(\bar\lambda\in[0,\lambda_{\min})\) such that for every \(i=1,\ldots,n\) and
		\((\theta,z)\in\R^p\times\mathcal Z\), $\|\nabla_\theta\ell(\theta,z)\|^2
		+\|z\|^4
		\leq
		C(\theta)(1+\Psi(z))$
		and $|\ell(\theta,z)|
		\leq
		L+\bar\lambda\|z-\xi_i\|^2.$
		
		\item For every compact set $\Xi \subset\mathcal Z$, there exist constants
		$a>0$ and $\bar{r}>0$ such that $\inf_{z\in \Xi}\pi_0\bigl(\closedball(z,r)\bigr)
		\ge a r^d$ for all $r\in(0,\bar{r}].$
%
		
	\end{enumerate}
\end{assumption}

The following theorem establishes a bound on the entropic approximation error over the parameter domain. Although the result is fairly standard, our formulation extends previous versions, such as \cite[Lemma 7]{liu2025stochastic}, to noncompact sample spaces. Under weak convexity, it further yields a quantitative graphical approximation of the subdifferential based on \cite{davis2022graphical}.

\begin{theorem}
\label{th:quantitative_approx_graph}
	Let $K=\Theta\times\Lambda$. Under \Cref{ass:WDRO_smooth_quadratic}, there exists $C>0$ such that, for all $\beta \in (0,e^{-1}]$,
\begin{enumerate}
    \item $ \sup_{(\theta,\lambda) \in \Theta \times \Lambda}
    |F^\beta(\theta,\lambda)-F(\theta,\lambda)|
    \leq C\beta |\log \beta |.$

    \item $\mathbb{H}_{\frac{1}{2\tau_\ell}}
    \left(
    \Graph_K (\nabla F^\beta+N_K),
    \Graph_K (\partial^c F+N_K) 
    \right)
    \leq
    \sqrt{\frac{2C \beta |\log \beta|}{\tau_\ell}},$
    where $\mathbb{H}_{\frac{1}{2\tau_\ell}}$ denotes the Hausdorff distance induced by the product norm $$(u,v) \mapsto \max \left\{\|u\|,\frac{1}{2\tau_\ell}\|v\|\right\}.$$
\end{enumerate}
\end{theorem}
\begin{proof}
    See \Cref{subsection:SGD_proofs}.
\end{proof}

\subsection{Correctness of the restricted problem}
\label{subsection:truncation}

We now discuss how to choose the interval \(\Lambda=[\lambda_{\min},\lambda_{\max}]\) so that restricting the joint minimization of the WDRO objective from \(\Theta\times[0,\infty)\) to \(\Theta\times\Lambda\) preserves its optimal value. We identify conditions under which this holds, depending on the growth of the loss relative to the transport cost. The lower bound on \(\lambda\) is governed by their relative asymptotic growth, while the upper bound follows from a simple coercivity argument in the dual variable. These conditions are primarily illustrative and are not necessarily
compatible with the growth conditions in
\Cref{ass:regularity_ell_c_merged}. Establishing a compatible restriction
may require stronger, problem-specific arguments.

\begin{proposition}[Lower bound on the dual solution]
	\label{prop:asymptotic_slope_lower_lambda}
	Let $F$ be defined by \eqref{eq:WDRO_objective_main}. For each \(\theta\in\Theta\), define the asymptotic slope
	\[
	\Gamma(\theta)
	:=
	\max_{1\le i\le n}\ 
	\limsup_{c(\xi_i,z)\to\infty}
	\frac{\ell(\theta,z)-\ell(\theta,\xi_i)}{c(\xi_i,z)}.
	\]
	Then every minimizer \((\theta^*,\lambda^*)\) of \(F\) on $\Theta \times \R_+$ satisfies $
	\lambda^*\ge \Gamma(\theta^*).$ In particular, if $\inf_{\theta \in \Theta} \Gamma(\theta)>0$, there exists $\lambda_{\min} >0$ such that $\min_{\Theta \times [0,\infty)} F = \min_{\Theta \times [\lambda_{\min},\infty)} F$.
\end{proposition}

\begin{proof}
	Fix \(\theta\in\Theta\) and suppose that \(\lambda<\Gamma(\theta)\). By the
	definition of \(\Gamma(\theta)\), there exist
	\(i\in\{1,\dots,n\}\), \(\varepsilon>0\), and a sequence
	\((z_k)_{k\in\mathbb N}\subset\mathcal Z\) such that
	\(c(\xi_i,z_k)\to\infty\) and, for all \(k\) large enough, $\ell(\theta,z_k)-\ell(\theta,\xi_i)
	\ge
	(\lambda+\varepsilon)c(\xi_i,z_k).$
	Therefore
	\[
	\ell(\theta,z_k)+\lambda\bigl(\rho-c(\xi_i,z_k)\bigr)
	\ge
	\ell(\theta,\xi_i)+\lambda\rho+\varepsilon c(\xi_i,z_k)
	\to+\infty.
	\]
	Hence $\sup_{z\in\mathcal Z}
	\Bigl\{\ell(\theta,z)+\lambda\bigl(\rho-c(\xi_i,z)\bigr)\Bigr\}
	=+\infty,$
	and consequently \(F(\theta,\lambda)=+\infty\). Thus no minimizer can satisfy
	\(\lambda^*<\Gamma(\theta^*)\). The second statement follows
	immediately if there is a positive uniform lower bound on $\Gamma(\theta)$, by taking $\lambda_{\min} = \inf_{\Theta} \Gamma /2$.
\end{proof}
The next proposition gives the complementary upper truncation.

\begin{proposition}[Upper bound on the dual variable]
	\label{prop:upper_truncation_lambda}
	Assume $c(z,z) = 0$ for all $z \in \mathcal{Z}$ and there exist \(\underline L\in\mathbb R\) and
	\((\bar\theta,\bar\lambda)\in\Theta\times[0,\infty)\) such that for all $\theta \in \Theta$, $\frac1n\sum_{i=1}^n \ell(\theta,\xi_i)\ge \underline L$ and $\overline V:=F(\bar\theta,\bar\lambda)<+\infty.$
	Then every minimizer \((\theta^*,\lambda^*)\) of \(F\) satisfies $\lambda^*\le
	\frac{\overline V-\underline L}{\rho}.$
\end{proposition}

\begin{proof}
	For every \((\theta,\lambda)\in\Theta\times[0,\infty)\), evaluating the
	inner supremum at \(z=\xi_i\) gives
	\[
	\sup_{z\in\mathcal Z}
	\Bigl\{\ell(\theta,z)+\lambda\bigl(\rho-c(\xi_i,z)\bigr)\Bigr\}
	\ge
	\ell(\theta,\xi_i)+\lambda\rho,
	\]
	because \(c(\xi_i,\xi_i)=0\). Therefore $F(\theta,\lambda)
	\ge
	\frac1n\sum_{i=1}^n \ell(\theta,\xi_i)+\lambda\rho
	\ge
	\underline L+\lambda\rho.$
	Let \((\theta^*,\lambda^*)\) be a minimizer. Since
	\(F(\theta^*,\lambda^*)\le F(\bar\theta,\bar\lambda)=\overline V\),
	we obtain $\underline L+\lambda^*\rho
	\le
	\overline V,$
	hence the result.
\end{proof}

\subsection{Stochastic gradient methods for regularized WDRO}
\label{subsection:gradient_methods}
We investigate stochastic gradient methods to minimize the WDRO objective \eqref{eq:WDRO_objective_main} based on the regularized objective \eqref{eq:regularized_objective}. For $m \geq 1$, an i.i.d sample $z_1, \ldots, z_m \sim \pi_0$ and $\beta >0$, we define the sampling approximation 
\begin{equation}
	\label{eq:H_beta_m}
	H_i^{\beta,m}(\theta,\lambda) = \beta \log \left( \frac{1}{m}\sum_{j=1}^m \exp  \left(
	\frac{\ell(\theta,z_j)+\lambda(\rho-c(\xi_i,z_j))}{\beta}
	\right) \right)
\end{equation}
and we set $F^{\beta,m}(\theta,\lambda) = \frac{1}{n} \sum_{i=1}^n H_i^{\beta,m}(\theta,\lambda).$

\paragraph{Fixed regularization and qualitative guarantees} First, we have qualitative guarantees for a stochastic gradient method with fixed regularization.

\begin{algorithm}[H]
\caption{Projected stochastic gradient method with fixed regularization}
\label{alg:SGD_fixed_reg}
\begin{algorithmic}[0]
\STATE \textbf{Generate} i.i.d. samples $z_1,\ldots,z_m\sim \pi_0$.
\STATE \textbf{Initialize} $(\theta_0,\lambda_0)\in \Theta\times\Lambda$ and positive stepsizes $(\alpha_k)_{k\in\mathbb{N}}$.
\FOR{$k\in\mathbb{N}$}
    \STATE Sample $i_k\sim \operatorname{Unif}\{1,\ldots,n\}$.
    \STATE
    $\begin{pmatrix}
        \theta_{k+1} \\
        \lambda_{k+1}
    \end{pmatrix}
    =
    \Pi_{\Theta\times\Lambda}
    \left(
    \begin{pmatrix}
        \theta_k \\
        \lambda_k
    \end{pmatrix}
    -
    \alpha_k
    \nabla H_{i_k}^{\beta,m}(\theta_k,\lambda_k)
    \right).$
\ENDFOR
\end{algorithmic}
\end{algorithm}
\begin{proposition} \label{prop:convergence_fixed_reg} Assume $F^{\beta,m}(\crit_{\Theta \times \Lambda} F^{\beta,m})$ has empty interior for any $\beta > 0$ and $m \in \N$. Under \Cref{ass:regularity_ell_c_merged}, let $\varepsilon > 0$. Then there exists $\bar{\beta} >0$ such that for any $\beta \in (0, \bar{\beta}]$, there exists $\Bar{m} \in \N$ such that for all $m \geq \bar{m}$, if $(\theta_k,\lambda_k)_{k \in \N}$ is generated by \Cref{alg:SGD_fixed_reg} with a sample size of $m$ and with stepsizes $\alpha_k > 0$, $\sum_{k=0}^\infty \alpha_k = \infty$ and $\sum_{k=0}^\infty \alpha_k^2 < \infty$, then  $$\limsup_{k\to\infty} \dist((\theta_k,\lambda_k), \crit_{\Theta \times \Lambda} F) \leq \varepsilon.$$
\end{proposition}

\begin{proof}
    This is a consequence of the approximation result  \Cref{th:gradient_conv_ell_c} and  standard convergence results for projected stochastic gradient methods. See, for example, \cite[Theorem 6.2, Chapter 5]{kushner2003stochastic}, \cite[Chapter 5.5]{borkar}, and also \cite[Theorem 5.1]{ermol1998stochastic}.
\end{proof}
The assumption that $F^{\beta,m}(\crit F^{\beta,m}_{\Theta \times \Lambda})$ has empty interior is common in stochastic approximation. Since the sampling approximation involves finite sums, this can be ensured for most functions used in practice, by the Morse--Sard theorem \cite{morse1939behavior,sard1942measure} and its generalizations \cite{shikhman,bolte2006nonsmooth}. 

\paragraph{Vanishing regularization schemes and convergence rates.}

Under \Cref{ass:WDRO_smooth_quadratic} there is a $\pi_0$-integrable function $a$ $\sup_{\theta\in\Theta}|\ell(\theta,z)|\leq a(z)$.
Fix an error tolerance $\bar b>0$ and deterministic regularization levels
$\beta_k>0$. By \Cref{lem:LLN_H_beta_m}, choose deterministic sample sizes
$m_k\geq m_{\beta_k,\bar b}$ before sampling, so that
\begin{equation}
\label{eq:m_k_bound}
\mathbb E\!\left[
\max_{1\leq i\leq n}\sup_{w\in K}
|H_i^{\beta_k,m_k}(w)-H_i^{\beta_k}(w)|
\,\middle|\,\mathcal F_k\right]\leq\bar b,
\qquad K=\Theta\times\Lambda,
\end{equation}
where $\mathcal F_k$ is the history before sampling at iteration $k$ and
new batches are independent of this history and of $i_k$.
The sample sizes are theoretical, since we do not provide a practical
rule for computing $m_{\beta_k,\bar b}$.
Developing such criteria is beyond the scope of this work, which focuses on smoothing. With this in mind, we consider a stochastic gradient method with varying regularization, described as follows.

\begin{algorithm}[H]
\caption{Theoretical sampling scheme for varying entropic smoothing}
\label{alg:SGD_vanishing_reg}
\begin{algorithmic}[0]
\STATE \textbf{Initialize} $(\theta_0,\lambda_0)\in \Theta\times\Lambda$ and choose step sizes $(\alpha_k)_{k\in\mathbb{N}}$ and regularization levels $(\beta_k)_{k \in \N}$.
\FOR{$k=0, \ldots N-1$}
    \STATE Sample $i_k\sim \operatorname{Unif}\{1,\ldots,n\}$.
    \STATE\textbf{Generate} a fresh i.i.d. batch $z_1,\ldots,z_{m_k}\sim\pi_0$, with $m_k$ chosen in advance as in \eqref{eq:m_k_bound}.
    \STATE $\begin{pmatrix}
        \theta_{k+1} \\
        \lambda_{k+1}
    \end{pmatrix}
    =
    \Pi_{\Theta\times\Lambda}
    \left(
    \begin{pmatrix}
        \theta_k \\
        \lambda_k
    \end{pmatrix}
    -
    \alpha_k
    \nabla H_{i_k}^{\beta_k,m_k}(\theta_k,\lambda_k)
    \right).$
\ENDFOR
\end{algorithmic}
\end{algorithm}
The corresponding convergence rates are as follows.
\begin{proposition} \label{prop:convergence_SGD_vanish_reg}
Under \Cref{ass:WDRO_smooth_quadratic} and the integrable loss-envelope condition above, let
\((\theta_k,\lambda_k)_{k \geq 0}\) be generated by \Cref{alg:SGD_vanishing_reg} with
$\alpha_k = \frac{\alpha_0}{\sqrt{k+1}} > 0$ and  $\beta_k = \frac{\beta_0}{k+1} \in (0,e^{-1}].$
For \(N \geq 1\), let \(R_N \in \{0,\ldots,N-1\}\) be drawn according to $\mathbb{P}(R_N = k)
=
\frac{\alpha_k}{\sum_{j=0}^{N-1}\alpha_j}.$
Then the following convergence rates hold:
\begin{enumerate}
    \item For any \(\tau > \tau_\ell\), where \(\tau_\ell\) is the weak convexity constant of
    \(\ell\) given by \Cref{ass:WDRO_smooth_quadratic}.2,
    \begin{equation*}
        \mathbb{E}
        \left[
        \left\|
        \nabla \varphi_{1/\tau}(\theta_{R_N},\lambda_{R_N})
        \right\|^2
        \right]
        =
        O\left(
        \frac{\log N}{\sqrt{N}} + \bar b
        \right),
    \end{equation*}
    where $\varphi_{1/\tau}(w)
    =
    \min_{w' \in \Theta \times \Lambda}
    \left\{
    F(w') + \frac{\tau}{2}\|w-w'\|^2
    \right\},$
    and \(\bar b\) is given by \eqref{eq:m_k_bound}.

    \item If, in addition, \(\ell(\cdot,z)\) is convex for \(\pi_0\)-almost every \(z\), then
    \begin{equation*}
        \mathbb{E}
        \left[
        F(\theta_{R_N},\lambda_{R_N})
        \right]
        -
        \min_{\Theta \times \Lambda} F
        =
        O\left(
        \frac{\log N}{\sqrt{N}} + \bar b
        \right).
    \end{equation*}
\end{enumerate}
\end{proposition}

\section{Entropic smoothing for sums of maximum functions}
\label{section:convergence_entropic_general}

To establish the guarantees stated in \Cref{sec:main_results}, we consider functions $H$ defined as pointwise maxima  and for $\beta>0$, their associated entropic regularization $H^\beta$ defined by
\begin{equation}
    \label{eq:H_and_Hbeta}
    H(w) = \max_{z \in \mathcal{Z}} h(w,z), \qquad H^\beta(w) = \beta \log \E_{z \sim \pi_{0}}[e^{h(w,z)/\beta}]
\end{equation}
Under differentiability assumptions on $h$, we investigate the convergence of $\nabla H^{\beta}$ to the subdifferential \eqref{eq:clarke} $\partial^c H$. Convergence and compactness results for the Gibbs measures as $\beta \to 0$ are developed in \Cref{lem:mass_near_maximal_subset} and \Cref{lem:expectation_limit_gibbs} in order to prove the graphical convergence results,  \Cref{th:graphical_convergence_one_batch} and  \Cref{cor:mini_batch_graphical_convergence}.

\subsection{On graphs of set-valued maps}
Since our analysis involves graphical convergence of gradient mappings to
set-valued subdifferential mappings, we first recall some useful notions for
set-valued maps. We will use the notation $\R^q \rightrightarrows \R^q$ to denote a map from $\R^q$ to the subsets of $\R^q$.

The following definition formally characterizes the perturbation appearing in \eqref{eq:contribution_convergence} and in our subsequent approximation results.

\begin{definition}[Graph-closed set-valued maps] \label{def:graph_closed} Let $\Theta \subset \R^q$ and let  $D : \Theta \rightrightarrows \R^q$ be a set-valued map. $D$ is called \emph{graph-closed} if its graph, defined as \eqref{eq:graph}, is a closed set. 
\end{definition}
Note that the Clarke subdifferential \eqref{eq:clarke} is graph-closed \cite{clarke1990optimization}. For $\delta > 0$, the $\delta$-enlargement of $D$ is the set-valued map $D^\delta$ where for each $w \in \R^q$
\begin{equation}
	\label{eq:graph_enlargment}
	D^{\delta}(w)
	:= \left\{\, v \in \mathbb{R}^q \ \middle|\ 
	\begin{array}{l}
		\exists (w',v') \text{ with } v' \in D(w'),\\
		\|v-v'\| \le \delta,\ \|w-w'\| \le \delta
	\end{array}
	\right\}.
\end{equation}
$D^\delta$ also corresponds to the set-valued map such that $\Graph D^{\delta} = \Graph D + (\delta \closedball) \times (\delta \closedball)$. Furthermore, if $D,D'$ are two set-valued maps with $D'\subset D^{\delta}$ for some $\delta >0$, then $\D(\Graph  D' , \Graph D) \leq \sqrt{2} \delta$, where $\D$ is the excess distance \eqref{eq:excess_distance}.

The following lemma shows that graphical approximation implies approximation of critical sets. This will be useful for our results on the convergence of critical sets.

\begin{lemma} \label{lem:graphconv_to_critconv} Let  $K \subset \R^q$ be a compact convex subset, let $D : \R^q \rightrightarrows \R^q$ be graph-closed and  locally bounded, and fix $\varepsilon > 0$. There exists $\delta > 0$ such that $\{w \in K \ : 0 \in D^{\delta}(w) + N_{K}(w) \} \subset \{w \in K \ : \ 0 \in D(w) + N_{K}(w)\} + \varepsilon \closedball$
\end{lemma}

\begin{proof} Assume toward a contradiction that there exists $\varepsilon > 0$ and  a positive sequence $(\delta_k)_{k \in \N}$ converging to $0$ such that for each $k \in \N$ there exists $w_k \in K$ satisfying $0 \in D^{\delta_k}(w_k) + N_{K}(w_k)$ and $\dist \left( w_k, \{w \in K \right.$ $\left. : \ 0 \in D(w) + N_{K}(w)\} \right) $ $> \varepsilon$. In particular, for each $k \in \N$, there exists $v_k \in D^{\delta_k}(w_k) \cap (- N_{K}(w_k))$. By compactness of $K$ and since $D$ is locally bounded, we may assume without loss of generality that $(w_k,v_k)$ converges to a pair $(\bar{w}, \bar{v}) \in \R^q \times \R^q$. For each $k$, take $(w_k', v_k') \in \Graph D$ such that $\|w_k - w_k'\| \leq \delta_k$ and $\|v_k - v_k'\| \leq \delta_k$. Then this new sequence also satisfies $(w_k', v_k') \to (\bar{w}, \bar{v})$. Hence since $D$ is graph-closed, $\bar{v} \in D(\bar{w})$. Moreover, since $-v_k \in N_K(w_k)$ for every $k$ and $N_K$ has a
closed graph, we have $-\bar v \in N_K(\bar w)$. Together with
$\bar v \in D(\bar w)$, this gives $0 \in D(\bar w)+N_K(\bar w)$.
Thus $\bar w$ belongs to the critical set
$S := \{w \in K : 0 \in D(w)+N_K(w)\}$. Consequently, $\dist(w_k,S) \leq \|w_k-\bar w\| \longrightarrow 0,$
contradicting $\dist(w_k,S)>\varepsilon$ for every $k$.
\end{proof}



\subsection{Gradient convergence for entropic regularization}

We now turn to the convergence result for entropic smoothing. We define the parameterized Gibbs measure by
\begin{equation}
    \label{eq:gibbs}
\di \pi_{\beta,w}(z) = \frac{e^{h(w,z)/\beta}}{Z_\beta(w)}\, \di \pi_0(z),
\qquad 
Z_\beta(w) = \int_{\mathcal{Z}} e^{h(w,z)/\beta}\, \di \pi_0( z).
\end{equation}
For the qualitative and general results, we work under the following assumptions. 

\begin{assumption}
\label{ass:regularity_h}
\begin{enumerate}
    \item[]
    \item $\supp(\pi_0)=\mathcal Z$;
    \item $(w,z)\mapsto h(w,z)$ is continuous on $\R^p\times \mathcal Z$;
    \item There exist a continuous coercive and $\pi_0$-integrable function $\psi:\mathcal Z\to \R$ and a locally bounded function $\kappa : \R^p \to [1, \infty)$ such that $h(w,z)\le -\kappa(w)\psi(z)$ for all $(w,z)\in \R^p\times \mathcal Z.$
\end{enumerate}
\end{assumption}

\begin{lemma}[Mass of a near maximizer set] \label{lem:mass_near_maximal_subset}
Under \Cref{ass:regularity_h}, for $w \in \R^p$ and $\varepsilon>0$, let $S(w,\varepsilon) := \{z \in \mathcal{Z} \ : \ h(w,z) > H(w) - \varepsilon \},$ where $H$ is defined as \eqref{eq:H_and_Hbeta}.  Then for any compact set $K \subset \R^p$, 
\[
\inf_{w \in K} \pi_{\beta,w}(S(w,2\varepsilon)) \geq 1 -  m_{K,\varepsilon}^{-1} e^{-\frac{\varepsilon}{\beta}},
\]
where $m_{K,\varepsilon}:=\inf_{w\in K}\pi_0\bigl(S(w,\varepsilon)\bigr)>0$ and $\pi_{\beta,w}$ is given by \eqref{eq:gibbs}.
\end{lemma}

\begin{proof} Note that the maximum function $H$ is well defined thanks to coercivity of $\psi$ given by \Cref{ass:regularity_h}.3. Fix $\varepsilon>0$ and a compact set $K\subset \R^p$. We first show that
\[
m_{K,\varepsilon}:=\inf_{w\in K}\pi_0\bigl(S(w,\varepsilon)\bigr)>0.
\]
For each $w\in K$, the set $S(w,\varepsilon)$ is open and nonempty by continuity of $h$ and the definition of $H(w)$. Since $\supp(\pi_0)=\mathcal Z$, every nonempty open subset of $\mathcal Z$ has strictly positive $\pi_0$-measure, so $\pi_0(S(w,\varepsilon))>0$ for all $w\in K$. Moreover, $w\mapsto \pi_0(S(w,\varepsilon))$ is lower semicontinuous on $K$. Indeed, if $w_n\to w$ and $z\in S(w,\varepsilon)$, then $h(w,z)>H(w)-\varepsilon$, and by continuity of $h$ and $H$, this strict inequality persists for all sufficiently large $n$, so eventually $z\in S(w_n,\varepsilon)$. Hence $\liminf_{n\to\infty}\mathbf 1_{S(w_n,\varepsilon)}(z)\ge \mathbf 1_{S(w,\varepsilon)}(z)$ for all $z\in\mathcal Z.$
Fatou's lemma then gives $\liminf_{n\to\infty}\pi_0\bigl(S(w_n,\varepsilon)\bigr)\ge \pi_0\bigl(S(w,\varepsilon)\bigr),$
hence $w\mapsto \pi_0(S(w,\varepsilon))$ is lower semicontinuous. Since it is strictly positive on the compact set $K$, we obtain $m_{K,\varepsilon}>0$. Now fix $w\in K$. We have
\[
\pi_{\beta,w}\bigl(S(w,2\varepsilon)^c\bigr)
=\frac{\int_{S(w,2\varepsilon)^c} e^{h(w,z)/\beta}\,\di\pi_0(z)}
{\int_{\mathcal Z} e^{h(w,z)/\beta}\,\di\pi_0(z)}
\le
\frac{\int_{S(w,2\varepsilon)^c} e^{h(w,z)/\beta}\,\di\pi_0(z)}
{\int_{S(w,\varepsilon)} e^{h(w,z)/\beta}\,\di\pi_0(z)}.
\]
For $z\in S(w,2\varepsilon)^c$, one has $h(w,z)\le H(w)-2\varepsilon$, while for $z\in S(w,\varepsilon)$, one has $h(w,z)>H(w)-\varepsilon$. Therefore
\[
\int_{S(w,2\varepsilon)^c} e^{h(w,z)/\beta}\,\di\pi_0(z)
\le e^{(H(w)-2\varepsilon)/\beta},
\]
and
\[
\int_{S(w,\varepsilon)} e^{h(w,z)/\beta}\,\di\pi_0(z)
\ge e^{(H(w)-\varepsilon)/\beta}\,\pi_0\bigl(S(w,\varepsilon)\bigr).
\]
Dividing the two inequalities yields
\[
\pi_{\beta,w}\bigl(S(w,2\varepsilon)^c\bigr)
\le \frac{e^{(H(w)-2\varepsilon)/\beta}}
{e^{(H(w)-\varepsilon)/\beta}\,\pi_0(S(w,\varepsilon))}
= \frac{1}{\pi_0(S(w,\varepsilon))}e^{-\varepsilon/\beta}
\le m_{K,\varepsilon}^{-1}e^{-\varepsilon/\beta}.
\]
Hence $\pi_{\beta,w}\bigl(S(w,2\varepsilon)\bigr)\ge 1-m_{K,\varepsilon}^{-1}e^{-\varepsilon/\beta}$ for all $w\in K,$ which proves the claim.
\end{proof}

\begin{lemma}[Uniform integrability] \label{lem:uniform_integrability_h}
Let $K\subset \R^p$ be compact. Under \Cref{ass:regularity_h}, there exists a convex non-decreasing function $\Phi : \R_+ \to \R$ with $\Phi(t)/t \xrightarrow[t \to \infty]{} \infty$, such that $$\sup_{\beta>0,\;w\in K}\int_{\mathcal Z}\Phi(|  \psi|)\,\di \pi_{\beta,w}<\infty.$$
In particular, the family of Gibbs measures $\{\pi_{\beta,w}: \beta>0,\ w \in K\}$ is uniformly integrable with respect to $\psi$, namely $\sup_{\beta >0, w \in K} \int_{\mathcal{Z}} |\psi| \mathbf{1}_{|\psi| > M} \di \pi_{\beta,w} \xrightarrow[M \to \infty]{} 0$,  and relatively compact.
\end{lemma}
\begin{proof}
Fix $\varepsilon>0$. For $w\in K$ let $S(w,2\varepsilon):=\{z\in\mathcal Z:\ h(w,z)>H(w)-2\varepsilon\}.$
De la Vall\'{e}e Poussin criterion \cite[Th.~22, p.~19]{dellacherie2011probabilities} gives us a convex function $\Phi : \R_+ \to \R$ such that $\Phi(t)/t \xrightarrow[t \to \infty]{} \infty$ and $\int_{\mathcal{Z}} \Phi (|\psi|)(z) \di \pi_0 < \infty$. We split
\[
\int_{\mathcal Z} \Phi(|\psi|)\,\di\pi_{\beta,w}
=
\int_{S(w,2\varepsilon)} \Phi(|\psi|)\,\di\pi_{\beta,w}
+
\int_{S(w,2\varepsilon)^c} \Phi(|\psi|)\,\di\pi_{\beta,w}.
\]
If $z\in S(w,2\varepsilon)$ then $H(w)-2\varepsilon<h(w,z)\le -\psi(z)$, hence
$\psi(z)<2\varepsilon-H(w)$. Since $H$ is continuous on the compact set $K$, letting
$\underline H_K:=\min_{w\in K}H(w)$ gives $\psi(z)\le 2\varepsilon-\underline H_K$. Since $\psi$ is also bounded below by coercivity, and $\Phi$ nondecreasing, there exists $C_{K,\varepsilon}<  \infty$ such that $\int_{S(w,2\varepsilon)}\Phi(|\psi|)\,\di\pi_{\beta,w}
\le C_{K,\varepsilon}.$
As to the integral over the complement, we write
\[
\int_{S(w,2\varepsilon)^c}\Phi(|\psi|)\,\di\pi_{\beta,w}
=
\frac{\int_{S(w,2\varepsilon)^c}\Phi(|\psi|)e^{h(w,z)/\beta}\,\di\pi_0(z)}
{\int_{\mathcal Z}e^{h(w,u)/\beta}\,\di\pi_0(u)}.
\]
We bound the numerator above, and the denominator below. For $z\in S(w,2\varepsilon)^c$ one has $h(w,z)\le H(w)-2\varepsilon$, hence
\[
\int_{S(w,2\varepsilon)^c}\Phi(|\psi(z)|) e^{h(w,z)/\beta}\,\di\pi_0(z)
\le
e^{(H(w)-2\varepsilon)/\beta}\int_{\mathcal Z}\Phi(|\psi(z)|)\,\di\pi_0(z).
\]
On the other hand, since $h(w,z)>H(w)-\varepsilon$ for $z\in S(w,\varepsilon)$,
\[
\int_{\mathcal Z}e^{h(w,u)/\beta}\,\di\pi_0(u)
\ge
\int_{S(w,\varepsilon)}e^{h(w,u)/\beta}\,\di\pi_0(u)
\ge
e^{(H(w)-\varepsilon)/\beta} m_{K,\varepsilon}
\]
where the lower bound constant $m_{K,\varepsilon}>0$ is given by \Cref{lem:mass_near_maximal_subset}. Dividing the two estimates gives
\[
\int_{S(w,2\varepsilon)^c}\Phi(|\psi|)\,\di\pi_{\beta,w}
\le
e^{-\varepsilon/\beta}
\frac{\int_{\mathcal Z}\Phi(|\psi|)\,\di\pi_0}{\pi_0(S(w,\varepsilon))}
\le
e^{-\varepsilon/\beta}
\frac{\int_{\mathcal Z}\Phi(|\psi|)\,\di\pi_0}{m_{K,\varepsilon}}.
\]
This yields the desired bound $\sup_{\beta>0,\;w\in K}
\int_{\mathcal Z}\Phi(|\psi|)\,\di\pi_{\beta,w}(z)
<
\infty.$ By De la Vall\'{e}e Poussin criterion \cite[Th.~22, p.~19]{dellacherie2011probabilities}, the family
$\{\pi_{\beta,w}:\beta>0,\ w\in K\}$ is uniformly integrable with respect to $\psi$. In particular, for any $M > 0$, $M \pi_{\beta,w}(|\psi|>M) \leq \int_{\mathcal{Z}}|\psi| \mathbf{1}_{|\psi| > M}  \di \pi_{\beta,w}$ hence we have $\sup_{\beta>0,w\in K}\pi_{\beta,w}(|\psi|>M)\xrightarrow[M \to \infty]{}0$ where $|\psi|$ is coercive,
which implies the family $(\pi_{\beta,w})_{\beta > 0, w \in K}$ is relatively compact, by Prokhorov's theorem.
\end{proof}
The following lemma gives a corresponding bound for the empirical measures.
\begin{lemma}[Uniform empirical moment bound]
	\label{lem:uniform_moment_empirical_gibbs}
	Under \Cref{ass:regularity_h}, let  $\varepsilon>0$ and $S(w,\varepsilon)
	:=
	\left\{
	z\in\mathcal Z:
	h(w,z)\geq H(w)-\varepsilon
	\right\}.$ Let $g:\mathcal Z\to[0,\infty)$ be continuous and suppose that $I_g
	:=
	\int_{\mathcal Z}g\,\di\pi_0<\infty.$
	Assume there exists $\varepsilon>0$ such that $p_{K,\varepsilon}
	:=
	\inf_{w\in K}
	\pi_0\bigl(S(w,\varepsilon)\bigr)>0$ and $R_{K,\varepsilon}
	:=
	 \sup_{w\in K, z\in S(w,2\varepsilon)}
	g(z)<\infty.$
	Then
	\[
	\sup_{w \in K,\ \beta>0,\ m\geq1}
	\mathbb E\left[
	\int_{\mathcal Z}g(z)\,
	\di\pi_{\beta,w}^{m}(z)
	\right]
	\leq
	R_{K,\varepsilon}
	+\frac{3I_g}{p_{K,\varepsilon}}
	<\infty,
	\]
	where $\pi_{\beta,w}^{m} \propto \sum_{j=1}^m e^{h(w,z_j)/\beta} \delta_{z_j}$, with $z_1,\ldots,z_m \overset{\mathrm{i.i.d.}}{\sim} \pi_0$.
\end{lemma}

\begin{proof}
	Fix $w\in K$, $\beta>0$, and $m\geq1$. We will use the shorthands
	\[
	A:=S(w,\varepsilon),
	\qquad
	B:=S(w,2\varepsilon),
	\qquad
	p:=\pi_0(A)\geq p_{K,\varepsilon}.
	\]
	For  $z_1,\ldots,z_m \overset{\mathrm{i.i.d.}}{\sim} \pi_0$ we define the number of samples in $S(w,\varepsilon)$ as $N_m:=\sum_{j=1}^m\mathbf 1_A(z_j).$ On the event $\{N_m\geq1\}$, the denominator of the empirical Gibbs
	measure  $\pi_{\beta,w}^m$ has the following lower bound.
	\begin{equation}
		\label{eq:moment_bound_denominator_gibbs}
		\sum_{j=1}^m e^{h(w,z_j)/\beta}
		\geq
		N_m e^{(H(w)-\varepsilon)/\beta}.
	\end{equation}
	Moreover, if a sample $z_j\notin B$, then
	\begin{equation}
		\label{eq:moment_bound_numerator_gibbs}
		e^{h(w,z_j)/\beta}
		\leq e^{(H(w)-2\varepsilon)/\beta}.
	\end{equation}
	Since $g\leq R_{K,\varepsilon}$ on $B$, by using \eqref{eq:moment_bound_denominator_gibbs} and \eqref{eq:moment_bound_numerator_gibbs}, we have on the event 	$\{N_m\geq1\}$,
	\begin{align}
		\int_{\mathcal Z} g(z)\,\di\pi_{\beta,w}^{m}(z) 
		& \leq R_{K,\varepsilon} + \frac{\sum_{j=1}^m e^{h(w,z_j)/\beta} g(z_j) \mathbf{1}_{B^c}(z_j)}{\sum_{j=1}^m e^{h(w,z_j)/\beta}}\\
		&\leq
		R_{K,\varepsilon}
		+
		\frac{e^{-\varepsilon/\beta}}{N_m}
		\sum_{j=1}^m g(z_j)\mathbf 1_{B^c}(z_j)
		\nonumber\\
		&\leq
		R_{K,\varepsilon}
		+
		\frac{1}{N_m}
		\sum_{j=1}^m g(z_j)\mathbf 1_{B^c}(z_j).
		\label{eq:empirical_gibbs_moment_on_N_positive}
	\end{align}
	For each $j$, define $N_{m,-j}:=
	\sum_{\substack{r=1\\r\neq j}}^m\mathbf 1_A(z_r).$
	Since \(A\subset B\) and $B^c \subset A^c$ , on the event \(\{z_j\in B^c\}\) one has
	\(N_m=N_{m,-j}\). Hence, since $N_{m,-j}$ and $z_j$ are independent, we easily have the inequalities
	\begin{align*}
		\mathbb E\left[
		\mathbf 1_{\{N_m\geq1\}}
		\frac{1}{N_m}
		\sum_{j=1}^m g(z_j)\mathbf 1_{B^c}(z_j)
		\right] & \leq
		mI_g\,
		\mathbb E\left[
		\frac{\mathbf 1_{\{N_{m,-1}\geq1\}}}{N_{m,-1}}
		\right] \\ & \leq
		2m I_g \,
		\mathbb E\left[
		\frac{1}{N_{m,-1}+1}
		\right].
	\end{align*}

	Now \(N_{m,-1}\) follows a binomial distribution with parameters $m-1$ and $p$. Since $N_{m,-1}\sim\operatorname{Bin}(m-1,p)$ and for $r \in \N$, $(r+1)^{-1}=\int_0^1t^r\,\mathrm dt$, we have
	\[
	\mathbb E\left[\frac{1}{N_{m,-1}+1}\right]
	=
	\int_0^1(1-p+pt)^{m-1}\,\mathrm dt
	=
	\frac{1-(1-p)^m}{mp}
	\leq\frac{1}{mp},
	\]
	where we used the probability generating function formula $\E[t^{N_{m,-1}}] = (1 - p + pt)^{m-1}$. Taking expectations in
	\eqref{eq:empirical_gibbs_moment_on_N_positive} gives
	\[
	\mathbb E\left[
	\mathbf 1_{\{N_m\geq1\}}
	\int_{\mathcal Z}g(z)\,
	\di\pi_{\beta,w}^{m}(z)
	\right]
	\leq
	R_{K,\varepsilon}
	+\frac{2I_g}{p}.
	\]
	On the event $\{N_m=0\}$, we have $\int_{\mathcal Z} g(z)\,\di\pi_{\beta,w}^{m}(z)
	\leq
	\sum_{j=1}^m g(z_j).$ Thus
	\begin{align*}
		\mathbb E\left[
		\mathbf 1_{\{N_m=0\}}
		\int_{\mathcal Z} g(z)\,
		\di\pi_{\beta,w}^{m}(z)
		\right]
		&\leq
		m I_g(1-p)^{m-1} \leq \frac{I_g}{p},
	\end{align*}
	where we can see $mp(1-p)^{m-1}\leq 1$  as a bound on the probability of exactly one success in a binomial experiment. Combining the estimates on $\{N_m\geq1\}$ and $\{N_m=0\}$ yields the desired result.
\end{proof}
\begin{lemma}[Empirical moment bound under growth assumptions]
	\label{lem:sufficient_uniform_moment_empirical_gibbs}
	Let $K\subset\mathbb R^p$ be compact and suppose that
	\Cref{ass:regularity_h} holds. Assume that there exist $q,c_K,r_K,a_K,\bar r>0$ such that if $M(w):=\argmax_{z\in\mathcal Z}h(w,z)$,
	\[
	h(w,z)
	\geq
	H(w)-c_K\dist(z,M(w))^q
	\quad\text{whenever}\quad
	\dist(z,M(w))\leq r_K,
	\]
     and $\pi_0\bigl(\closedball(m,r)\bigr)
	\geq a_Kr^d$ for all $w\in K,\ m\in M(w),\
	r\in(0,\bar r].$
	Let $g:\mathcal Z\to[0,\infty)$ be continuous and satisfy $I_g:=\int_{\mathcal Z}g(z)\,\di\pi_0(z)<\infty.$
	Then, for all $\varepsilon>0$,
	\[
	\sup_{w\in K,\ \beta>0,\ m\geq1}
	\mathbb E\left[
	\int_{\mathcal Z}g(z)\,
	\di\pi_{\beta,w}^m(z)
	\right]
	\leq
	R_{K,\varepsilon}
	+\frac{3I_g}{a_Kr_\varepsilon^d}
	<\infty,
	\]
	where $r_\varepsilon:=\min\{r_K,\bar r,
	(\varepsilon/c_K)^{1/q}\}$,
	$\underline H_K:=\min_{w\in K}H(w)$, and
	$R_{K,\varepsilon}:=
	\sup\{g:z\in \{\psi\leq\max\{0,2\varepsilon-\underline H_K\}\}\}$.
\end{lemma}

\begin{proof}
	Since $H$ is the pointwise supremum of the continuous functions
	$h(\cdot,z)$, it is lower semicontinuous. Hence $M(w)$ is nonempty since $K$ is compact and  $\underline H_K>-\infty$. Fix $\varepsilon>0$, $w\in K$, and $m_w\in M(w)$. We set $S(w,\varepsilon)
	:=
	\left\{
	z\in\mathcal Z:
	h(w,z)\geq H(w)-\varepsilon
	\right\}.$ By the growth
	condition on $h$, $\closedball(m_w,r_\varepsilon)
	\subset S(w,\varepsilon)$. Hence the condition on $\pi_0$ gives us $\pi_0\bigl(S(w,\varepsilon)\bigr)
	\geq a_Kr_\varepsilon^d,$
	and consequently $p_{K,\varepsilon}
	\geq a_Kr_\varepsilon^d>0.$
	Now, if $z\in S(w,2\varepsilon)$, then
	\[
	\underline H_K-2\varepsilon
	\leq
	H(w)-2\varepsilon
	\leq
	h(w,z)
	\leq
	-\kappa(w)\psi(z),
	\]
	so that $\kappa(w)\psi(z)
	\leq
	2\varepsilon-\underline H_K.$
	Since $\kappa(w)\geq1$, this implies $\psi(z)
	\leq
	\max\{0,2\varepsilon-\underline H_K\}.$
	Thus $S(w,2\varepsilon)
	\subset
	\left\{
	z\in\mathcal Z:
	\psi(z)\leq
	\max\{0,2\varepsilon-\underline H_K\}
	\right\}.$
	The set on the right is compact by coercivity of $\psi$, and hence
	$R_{K,\varepsilon}<\infty$ by continuity of $g$. The conclusion follows
	from \Cref{lem:uniform_moment_empirical_gibbs}.
\end{proof}

\begin{lemma} \label{lem:gibbs_measure_cv_max}
Let $K\subset \R^p$ be compact. Under \Cref{ass:regularity_h}, let
$(\beta_k,w_k)_{k\in\N}$ be any sequence from $(0,\infty)\times K$
converging to $(0,w)\in [0,\infty)\times K$.
Then any weak accumulation point $\pi$ of
$(\pi_{\beta_k,w_k})_{k\in\N}$ is supported on
$\argmax_{z\in\mathcal Z} h(w,z)$.
\end{lemma}
\begin{proof}
The family $(\pi_{\beta_k,w_k})_{k \in \N}$ is relatively compact for weak convergence, hence without loss of generality, after passing to a subsequence (\Cref{lem:uniform_integrability_h}), we may assume that
$\pi_{\beta_k,w_k}\to\pi$ weakly. Fix $\varepsilon>0$ and define the closed near-maximal set
\[
\overline S(w,\varepsilon):=\{z\in\mathcal Z:\ h(w,z)\ge H(w)-\varepsilon\}.
\]
We first show that for all sufficiently large $k$,
$S(w_k,\varepsilon)\subset \overline S(w,2\varepsilon).$
If not, there exist a subsequence $(w_{k_n})_{n \in \N}$ and a sequence $(z_n)_{n \in \N} \subset \mathcal Z$ such that
$h(w_{k_n},z_n)>H(w_{k_n})-\varepsilon$ while $h(w,z_n)<H(w)-2\varepsilon$. By Assumption~\ref{ass:regularity_h}.3 we have
$h(w_{k_n},z_n)\le -\kappa(w_{k_n})\psi(z_n)\le -\psi(z_n)$.
Since $(w_k)_{k \in \N}\subset K$ and $H$ is continuous on $K$, the sequence $(H(w_k))_{k \in \N}$ is bounded below.
Hence $-\psi(z_n)\ge h(w_{k_n},z_n)>H(w_{k_n})-\varepsilon$, so $(\psi(z_n))$ is bounded.
Because $\psi$ is coercive, we may assume $z_n\to z$. Passing to the limit in $h(w_{k_n},z_n)>H(w_{k_n})-\varepsilon$ and using continuity of $h$ and $H$
gives $h(w,z)\ge H(w)-\varepsilon$.
On the other hand, $h(w,z_n)<H(w)-2\varepsilon$ for all $n$, and continuity of $h(w,\cdot)$ yields
$h(w,z)\le H(w)-2\varepsilon$ which is a contradiction.

By \Cref{lem:mass_near_maximal_subset}, there exists $C_{K,\varepsilon}>0$ such that $\pi_{\beta_k,w_k}(S(w_k,\varepsilon))
\ge 1-C_{K,\varepsilon}e^{-\varepsilon/(2\beta_k)}.$
Hence for all sufficiently large $k$, the inclusion above gives
\[
\pi_{\beta_k,w_k}(\overline S(w,2\varepsilon))
\ge
\pi_{\beta_k,w_k}(S(w_k,\varepsilon))
\ge
1-C_{K,\varepsilon}e^{-\varepsilon/(2\beta_k)}.
\]
In particular, since $\beta_k\to0$, the right-hand side converges to $1$. Because $\overline S(w,2\varepsilon)$ is closed, the Portmanteau theorem gives
\[
\pi(\overline S(w,2\varepsilon))
\ge
\limsup_{k\to\infty}\pi_{\beta_k,w_k}(\overline S(w,2\varepsilon))
=1.
\]
Hence $\pi(\overline S(w,2\varepsilon))=1$. Since this holds for every $\varepsilon>0$, we have
\[
\supp(\pi)\subset
\bigcap_{\varepsilon>0}\overline S(w,2\varepsilon)
=
\{z\in\mathcal Z:\ h(w,z)=H(w)\}
=
\argmax_{z\in\mathcal Z} h(w,z).
\]
\end{proof}

\begin{lemma}
\label{lem:expectation_limit_gibbs}
Under \Cref{ass:regularity_h}, let $K \subset \R^p$ be compact and  let
$g : \R^p \times \mathcal{Z} \to \R$ be jointly continuous. Assume there exists a locally bounded function
$C : \R^p \to [1,\infty)$ and a $\pi_0$-integrable and continuous coercive function $\psi : \mathcal{Z} \to \R_+$ such that $\|g(w,z)\| \le C(w)\psi(z)$ for all $(w,z)\in \R^p\times\mathcal{Z}.$ For $w\in\R^p$, let $M(w) := \argmax_{z\in\mathcal Z} h(w,z),$
and define
\begin{equation}
\label{eq:expectation_M}
\E_{M(w)}[g(w)]
:=
\left\{
\int_{\mathcal Z} g(w,z)\,dP(z)
:\;
P\in\mathcal P(\mathcal Z),\; P(M(w))=1
\right\}.
\end{equation}
Then $\D
\big(
\Graph_K \E_{\pi_{\beta,\cdot}}[g],
\Graph_K \E_{M(\cdot)}[g]
\big)
\longrightarrow 0$ as $\beta \to 0^+.$
\end{lemma}\begin{proof}
We treat the case where $g$ is real-valued; the vector-valued case follows by applying the same argument componentwise. Let $(\beta_k,w_k)_{k\in\N}\subset (0,\infty)\times K$ be such that $(\beta_k,w_k)\to(0,w)$ for some $w\in K$. Our goal is to show that
\[
\dist\!\Big(
\E_{z\sim \pi_{\beta_k,w_k}}[g(w_k,z)],
\E_{M(w)}[g(w,\cdot)]
\Big) \xrightarrow[k \to \infty]{} 0.
\]
where the expectation $\E_{M(w)}$ is defined in \eqref{eq:expectation_M}. Suppose, for contradiction, that there exist $\varepsilon>0$ and a subsequence, still denoted $(\beta_k,w_k)$, such that
\begin{equation}
    \label{eq:lower_bound_eps}
    \dist\!\Big(
\E_{z\sim \pi_{\beta_k,w_k}}[g(w_k,z)],
\E_{M(w)}[g(w,\cdot)]
\Big)\ge \varepsilon
\qquad \forall k\in\N.
\end{equation}
Set $\widetilde\pi_k:=\pi_{\beta_k,w_k}$. By \Cref{lem:uniform_integrability_h}, the family $(\widetilde\pi_k)_{k\in\N}$ is relatively compact, so, up to extraction, $\widetilde\pi_k\to\pi$ weakly  for some probability measure $\pi$. By \Cref{lem:gibbs_measure_cv_max}, the measure $\pi$ is supported on $M(w)$. Therefore it is enough to prove that
\begin{equation}
\label{eq:conv_key_g}
\int_{\mathcal Z} g(w_k,z)\,\di\widetilde\pi_k(z)
\xrightarrow[k \to \infty]{}
\int_{\mathcal Z} g(w,z)\,\di\pi(z),
\end{equation}
for then the right-hand side belongs to $\E_{M(w)}[g(w,\cdot)]$, contradicting \eqref{eq:lower_bound_eps}. Since $K$ is compact and $C$ is locally bounded, we may set $C_K:=\sup_{u\in K}C(u)<\infty,$
so that $|g(u,z)|\le C_K\psi(z)$ for all $(u,z)\in K\times\mathcal Z.$ Now $w_k\to w$ and $\widetilde\pi_k\to\pi$ weakly, hence Slutsky's theorem yields $\delta_{w_k}\otimes\widetilde\pi_k \to \delta_w\otimes\pi$ weakly on $\R^p\times\mathcal Z.$
For $R>0$, define the truncation $g_R(u,z):=\max(-R,\min(g(u,z),R)).$
Then $g_R$ is bounded and continuous on $\R^p\times\mathcal Z$, so
\[
\int_{\R^p\times\mathcal Z} g_R(u,z)\,\di(\delta_{w_k}\otimes\widetilde\pi_k)(u,z)
\longrightarrow
\int_{\R^p\times\mathcal Z} g_R(u,z)\,\di(\delta_w\otimes\pi)(u,z),
\]
that is, $\int_{\mathcal Z} g_R(w_k,z)\,\di\widetilde\pi_k(z)
\longrightarrow
\int_{\mathcal Z} g_R(w,z)\,\di\pi(z).$ On the other hand,
\[
|g(u,z)-g_R(u,z)|
\le |g(u,z)|\,\mathbf 1_{\{|g(u,z)|>R\}}
\le C_K\psi(z)\,\mathbf 1_{\{\psi(z)>R/C_K\}}.
\]
By \Cref{lem:uniform_integrability_h}, we have the uniform integrability $\sup_k\int_{\mathcal Z}\psi \mathbf{1}_{\{\psi > M\}}\,\di\widetilde\pi_k \xrightarrow[M \to \infty]{} 0.$ Therefore
\[
\sup_k\Big|
\int_{\mathcal Z}(g-g_R)(w_k,z)\,\di\widetilde\pi_k(z)
\Big|
\le
C_K\sup_k\int_{\mathcal Z}\psi\mathbf 1_{\{\psi>R/C_K\}}\,\di\widetilde\pi_k
\xrightarrow[R \to\infty]{} 0
\]
Also, since $\psi$ is nonnegative and continuous, the Portmanteau theorem gives
$\int_{\mathcal Z}\psi(z)\,\di\pi(z)
\le
\liminf_{k\to\infty}\int_{\mathcal Z}\psi(z)\,\di\widetilde\pi_k(z)
<\infty,$ and hence the same estimate yields $
\int_{\mathcal{Z}}|g(w,z)-g_R(w,z)|\,\di\pi(z)\longrightarrow 0$ as $R \to \infty$. For every $R>0$,
\[
\Big|
\int_{\mathcal Z} g(w_k,z)\,\di\widetilde\pi_k(z)
-
\int_{\mathcal Z} g(w,z)\,\di\pi(z)
\Big|
\le
\Big|\int_{\mathcal Z} (g-g_R)(w_k,z)\,\di\widetilde\pi_k(z)\Big|
\]
\[
\qquad\quad
+
\Big|\int_{\mathcal Z} g_R(w_k,z)\,\di\widetilde\pi_k(z)
-
\int_{\mathcal Z} g_R(w,z)\,\di\pi(z)\Big|
+
\Big|\int_{\mathcal Z} (g_R-g)(w,z)\,\di\pi(z)\Big|.
\]
Taking $\limsup_{k\to\infty}$ and then letting $R\to\infty$, we obtain \eqref{eq:conv_key_g}. This yields the desired contradiction, hence the result.
\end{proof}

Now we may derive the main result of this part. The result requires the following additional differentiability and integrability assumptions on $h$.

\begin{assumption}[Integrable gradient]
\label{ass:differentiable_h}
    \begin{enumerate}
        \item[]
        \item For all $z \in \mathcal{Z}$, $h(\cdot, z)$ is differentiable and $\nabla_w h$ is jointly continuous.   
        \item There exists a  $\pi_0$-integrable function $\varphi \ : \ \mathcal{Z} \to [0, \infty)$ and a locally bounded function $C \ : \ \R^p \to [1, \infty)$ such that $\int_\mathcal{Z} \varphi \di \pi_0 < \infty$ and  $\|\nabla_w h(w,z)\| \leq C(w) \varphi(z)$ for all $(w,z) \in \R^p \times \mathcal{Z}$.
    \end{enumerate}
\end{assumption}

\begin{theorem} \label{th:graphical_convergence_one_batch} Under \Cref{ass:regularity_h} and \Cref{ass:differentiable_h}, for any compact set $K \subset \R^p$, $$\D \left(\Graph_K \E_{z \sim \pi_{\beta, (\cdot)}}[\nabla_w h(\cdot,z)] , \Graph_K \partialc H \right) \xrightarrow[\beta \to 0^+]{} 0.$$
\end{theorem}
\begin{proof} Under \Cref{ass:regularity_h} and \Cref{ass:differentiable_h},   the envelope theorem  {\cite[Th. 2.8.2]{clarke1990optimization}} applies to give  $$\partialc H(w) = \left\{\E_{z \sim \pi}[\nabla_w h(w,z)] \ : \ \supp \pi \subset \argmax_{z \in \mathcal{Z} }h(w,z) \right\}.$$
The result is then a direct consequence of \Cref{lem:expectation_limit_gibbs} applied to $g = \nabla_w h$.
\end{proof}

\begin{corollary}[Mini-batch regularization] \label{cor:mini_batch_graphical_convergence} Let $h_1, \ldots, h_n : \R^p \times \mathcal{Z} \to \R$. Suppose Assumptions \ref{ass:regularity_h} and \ref{ass:differentiable_h} hold for each $h_i$, $i=1, \ldots,n$. Let $H_i = \max_{z \in \mathcal{Z}} h_i(\cdot, z)$,  $F := \frac{1}{n}\sum_{i=1}^n H_i$ and $F^\beta  =  \frac{1}{n}  \sum_{i = 1 }^n H_{i}^{\beta}$ where  each $H_i^\beta$ is defined as \eqref{eq:H_and_Hbeta}. Then for any compact set $K \subset \R^p$, $$\D \left( \Graph_K \nabla F^\beta, \Graph_K  \partialc F \right) \xrightarrow[\beta \to 0^+]{} 0$$
\end{corollary}
\begin{proof}
Let $K$ be a compact set. Suppose toward a  contradiction that  $$\mathbb{D}\!\left(\Graph_K \nabla  F^\beta,\Graph_K \partialc F\right)$$ does not converge to $0$. Then there exist $e>0$, $\beta_k\downarrow0$, and $(x_k,v_k)\in\Graph_K\nabla F^{\beta_k}$ such that $\dist\!\big((x_k,v_k),\Graph_K\partialc F\big)\ge e.$
Since $K$ is compact, we may assume $x_k\to\bar x\in K$. We set $v_k=\tfrac1n\sum_{i=1}^n\nabla H_{i}^{\beta_k}(x_k).$
\Cref{th:graphical_convergence_one_batch} gives the graphical convergence
$$\mathbb D(\Graph_K\nabla H_{i}^{\beta},\Graph_K\partialc H_i)\to0,$$ which yields $(y_{i,k},u_{i,k})\in\Graph_K\partialc H_i$ with $\|(x_k,\nabla H_{i}^{\beta_k}(x_k))-(y_{i,k},u_{i,k})\|\to0.$
Hence $y_{i,k}\to\bar x$ and $u_{i,k}-\nabla H_{i}^{\beta_k}(x_k)\to0$. By local boundedness and closedness of $\Graph\partialc H_i$, we may pass to a subsequence such that $u_{i,k}\to\bar u_i\in\partialc H_i(\bar x)$. Consequently $v_k\to\bar v:=\tfrac1n\sum_{i=1}^n\bar u_i\in\tfrac1n\sum_{i=1}^n\partialc H_i(\bar x).$ By the envelope theorem {\cite[Th. 2.8.2]{clarke1990optimization}}, each max function $H_i$ is \emph{Clarke regular} in the sense of \cite[Definition 2.3.4]{clarke1990optimization}. Hence the sum rule $\tfrac1n\sum_{i=1}^n\partialc H_i(\bar x) = \partialc F(\bar{x})$ holds and $(\bar x,\bar v)\in\Graph_K\partialc F$. Since $(x_k,v_k)\to(\bar x,\bar v)$, we obtain $\dist((x_k,v_k),\Graph_K\partialc F)\to0,$
which contradicts $$\dist((x_k,v_k),\Graph_K\partialc F)\ge e.$$
\end{proof}

\subsection{Quantitative analysis}
\label{subsection:quantitative_analysis}
Throughout this part we fix a compact set $K \subset \R^p$. We produce quantitative bounds on the function gap $H - H^\beta$  and the graphical distance between $\nabla H^{\beta}$ and $\partialc H$. Compared to the qualitative results, these complementary guarantees rely on stronger  regularity and geometric conditions which we state below. 
\begin{assumption}[Integrable gradient]
\label{ass:L_smooth_h}
    \begin{enumerate}
        \item[]
        \item $h$ is jointly continuous and for all $z \in \mathcal{Z}$, $h(\cdot, z)$ is differentiable.
        \item There exist a continuous $\pi_0$-integrable function $\Psi \ : \ \mathcal{Z} \to [0, \infty)$ and a locally bounded function $C \ : \ \R^p \to [0, \infty)$ such that $|h(w,z)|+ \|\nabla_w h(w,z)\| \leq C(w) \Psi(z)$ for all $(w,z) \in \R^p \times \mathcal{Z}$.
    \end{enumerate}
\end{assumption}

\begin{assumption}[Growth around maximizers] \label{ass:error_bound}
   There exist constants $q>0$, $c_K>0$, $r_K>0$, $a_K>0$ such that for every $w\in K$, 
   \begin{enumerate}
    \item  \(\mathcal Z\subset\R^d\), $M(w):=\argmax_{z\in\mathcal Z} h(w,z)$ is nonempty and for all $z\in\mathcal Z$ such that  $\dist(z,M(w))\le r_K$, $
h(w,z)\ge H(w)- c_K\,\dist(z,M(w))^q$.
\item  For every $m\in M(w)$ and every $r\in(0,r_K]$, $\pi_0(\closedball (m,r))\ge a_K r^d.$
\end{enumerate}
\end{assumption}

\begin{lemma}
	\label{lem:value_gap_one_batch}
	Under \Cref{ass:L_smooth_h} and \Cref{ass:error_bound}, there exists \(C_K>0\) such
	that, for every \(0<\beta\leq e^{-1}\),
	\[
	\sup_{w\in K}\delta_\beta(w)
	\leq C_K\beta|\log\beta|, 
	\qquad \text{where} \quad 
	\delta_\beta(w):=H(w)-H^\beta(w)\geq0.
	\]
\end{lemma}

\begin{proof}
	Fix \(w\in K\). Since \(h(w,z)\leq H(w)\) for every
	\(z\in\mathcal Z\), we have $H^\beta(w)
	=
	\beta\log
	\int_{\mathcal Z}
	e^{h(w,z)/\beta}\,\di\pi_0(z)
	\leq H(w).$ Therefore, \(\delta_\beta(w)\geq0\). For \(\varepsilon>0\), define the near-maximizer set as 
	\[
	A_\varepsilon(w)
	:=
	\left\{
	z\in\mathcal Z:
	h(w,z)\geq H(w)-\varepsilon
	\right\}.
	\]
	Restricting the integral defining \(H^\beta\) to
	\(A_\varepsilon(w)\), we obtain
	\begin{align*}
		H^\beta(w)
		&\geq
		\beta\log
		\int_{A_\varepsilon(w)}
		e^{h(w,z)/\beta}\,\di\pi_0(z)\\
		&\geq
		\beta\log
		\left(
		e^{(H(w)-\varepsilon)/\beta}
		\pi_0(A_\varepsilon(w))
		\right)\\
		&=
		H(w)-\varepsilon
		+\beta\log\pi_0(A_\varepsilon(w)).
	\end{align*}
	Equivalently,
	\begin{equation}
		\label{eq:value_gap_near_maximal_mass}
		\delta_\beta(w)
		\leq
		\varepsilon
		+
		\beta\log
		\left(
		\frac{1}{\pi_0(A_\varepsilon(w))}
		\right).
	\end{equation}
	Now, set $r_\varepsilon
	:=
	\min\left\{
	r_K,
	\left(\frac{\varepsilon}{c_K}\right)^{1/q}
	\right\}.$
	For any \(m\in M(w)\) and
	\(z\in\closedball(m,r_\varepsilon)\), we have $\dist(z,M(w))
	\leq \|z-m\|
	\leq r_\varepsilon.$
	Hence, by \Cref{ass:error_bound},
	\[
	h(w,z)
	\geq
	H(w)-c_K\dist(z,M(w))^q
	\geq
	H(w)-c_Kr_\varepsilon^q
	\geq
	H(w)-\varepsilon.
	\]
	This means $\closedball(m,r_\varepsilon)
	\subset A_\varepsilon(w).$
	\Cref{ass:error_bound}.2 therefore gives $\pi_0(A_\varepsilon(w))
	\geq
	a_Kr_\varepsilon^d.$
	Let \(\alpha:=d/q\). Taking \(\varepsilon=\beta\) in
	\eqref{eq:value_gap_near_maximal_mass}, we find
	\begin{align*}
		\delta_\beta(w)
		\leq
		\beta
		+
		\beta\log
		\left(
		\frac{1}{a_Kr_\beta^d}
		\right)
		=
		\beta
		+
		\beta\log
		\left[
		a_K^{-1}
		\max\left\{
		r_K^{-d},
		c_K^\alpha\beta^{-\alpha}
		\right\}
		\right].
	\end{align*}
	Define $C_0
	:=
	\max\left\{
	1,\,
	a_K^{-1}r_K^{-d},\,
	a_K^{-1}c_K^\alpha
	\right\}.$
	Since \(0<\beta\leq e^{-1}<1\), we have $\frac{1}{a_Kr_\beta^d} = a_K^{-1} \max \{r_K^{-d}, \beta^{- \alpha}  c_K^{\alpha}\}
	\leq
	C_0\beta^{-\alpha}.$
	Thus, $\delta_\beta(w)
	\leq
	\beta+\beta\log C_0+\alpha\beta|\log\beta|.$
	Moreover, \(|\log\beta|\geq1\) and $C_0 \geq 1$ so
	\[
	\delta_\beta(w)
	\leq
	\bigl(1+\log C_0+\alpha\bigr)
	\beta|\log\beta|.
	\]
	Taking the supremum over \(w\in K\) proves the result.
\end{proof}
Now, let us state the main result of this part.

\begin{proposition} \label{prop:quantitative_rate_h}  Let $h_1, \ldots, h_n : \R^p \times \mathcal{Z} \to \R$ and let $K \subset \R^p$ be a nonempty compact convex set. Suppose Assumptions \ref{ass:L_smooth_h} and \ref{ass:error_bound} hold for each $h_i$, $i=1, \ldots,n$. Assume there exists $\tau_h >0$ such that for all $z \in \mathcal{Z}$, each $h_i(\cdot,z)$ is $\tau_h$-weakly convex on $K$. Let $H_i = \max_{z \in \mathcal{Z}} h_i(\cdot, z)$,  $F := \frac{1}{n}\sum_{i=1}^n H_i$ and  $F^\beta  =  \frac{1}{n}  \sum_{i = 1 }^n H_{i}^{\beta}.$ Then there exists $C_K > 0$ such that for all $\beta \in (0, e^{-1}]$,
\begin{enumerate}
    \item $\sup_{w \in K} |F^\beta(w) - F(w)| \leq C_K \beta |\log \beta |$
    \item $\mathbb{H}_{\frac{1}{2\tau_h}} \left( \Graph_K (\nabla F^\beta+N_K), \Graph_K (\partialc F+N_K) \right)  \leq \sqrt{\frac{2C_K \beta |\log \beta|}{\tau_h}}$, where  $\mathbb{H}_{\frac{1}{2\tau_h}}$ is the Hausdorff distance induced by the product norm $(u,v) \mapsto \max \left\{\|u\|, \frac{1}{2\tau_h} \|v\| \right\}$.
\end{enumerate}
\end{proposition}
\begin{proof}
	 The first item is a direct consequence of \Cref{lem:value_gap_one_batch}. For the second item, apply \cite[Theorem 5.1]{davis2022graphical} to $F+\iota_K$ and $F^\beta+\iota_K$, where $\iota_K$ is the convex indicator of $K$. These functions have common domain $K$, and their subdifferentials are $\partialc F+N_K$ and $\nabla F^\beta+N_K$, respectively.
\end{proof}

\section{Approximate stochastic gradient methods}
\label{sec:convergence_rates_general}
In this section, we fix a compact convex set $K \subset \R^p$ and we consider an inexact stochastic gradient method for minimizing a finite sum of max functions $F = \frac{1}{n} \sum_{i=1}^n H_i$ on $K$, which corresponds to the setting of \Cref{prop:quantitative_rate_h}. We consider an inexact stochastic subgradient method initialized at $w_0 \in K$ where for each
iteration $k$, the index $i_k \in \{1,\ldots,n\}$ is sampled uniformly and
\begin{equation}
    \label{eq:algo_sgd_beta}
    w_{k+1}
    =
    \Pi_K\bigl(w_k-\alpha_k\nabla \widetilde H_{i_k}^{k}(w_k)\bigr),
\end{equation}
where for each
$i \in \{1,\ldots,n\}$ and $k\in\mathbb{N}$, $\widetilde H_i^k:K\to\mathbb{R}$ is a  differentiable approximation of $H_i^{\beta_k}$, given the current regularization level $\beta_k > 0$.   Under (weak) convexity assumptions, we derive bounds for general schedules and vanishing regularization coefficients  $(\beta_k)_{k \in \mathbb{N}}$. In particular, up to errors due to sampling approximations, we obtain a rate of $O(\log(N)/\sqrt{N})$, as $N \to \infty$. The proofs rely on the quantitative estimates established in \cref{subsection:quantitative_analysis}. For the weakly convex case, \Cref{prop:rates_weakly_cvx}, we adapt proofs from \cite{davis2019stochastic}. We set for all $k \geq 0$ and $\beta > 0$
\[
    \widetilde F_k(w)
    :=
    \frac{1}{n}\sum_{i=1}^n \widetilde H_i^k(w)
    \quad \text{and} \quad
    F^{\beta}(w)
    :=
    \frac{1}{n}\sum_{i=1}^n H_i^{\beta}(w).
\]
We work under the following set of assumptions.
\begin{assumption}[Approximation errors]
\label{ass:bounds_noise_beta}
\begin{enumerate}
    \item[]
    \item For all $k\in\mathbb{N}$, $i \in \{1,\ldots,n\}$, there exists a measurable function $b_{i,k} : \R^p \to \R_+$ such that for all $w\in K$, $$\mathbb{E}
        \left[
            |\widetilde H_i^k(w_k)-H_i^{\beta_k}(w_k)| | w_k
        \right]
        \leq b_{i,k}(w_k).$$
    \item There exists $\widetilde M_K^2\geq 0$ such that $\sup_{k\in\mathbb{N}}
        \mathbb{E}
        \left[
            \|\nabla \widetilde H_{i_k}^{k}(w_k)\|^2
            \mid w_k
        \right]
        \leq \widetilde M_K^2 .$
    \item There exists $C_K>0$ such that, for every $\beta>0$ sufficiently
    small, $$\sup_{w\in K}|F^\beta(w)-F(w)| \leq C_K\,\beta|\log\beta|.$$
\end{enumerate}
\end{assumption}
Conditions 1 and 2 capture the natural sampling noise associated with $\pi_0$. Condition 3 serves as a blanket assumption controlling the regularization error. \Cref{prop:quantitative_rate_h} establishes this bound in our WDRO setting. Under \Cref{ass:bounds_noise_beta} we define the model error by
\begin{equation}
	\label{eq:b_bar}
    \bar b_k(w)
:=
\frac{1}{n}\sum_{i=1}^n b_{i,k}(w),
\qquad \text{for all }w\in K.
\end{equation}

\subsection{Convex case}

In this part, the stochastic models
\(\widetilde H_i^k\) are assumed to be convex. This assumption is motivated by
the observation that, in practice, sampling-based approximations of
\(H_i^{\beta_k}\), such as \eqref{eq:regularized_objective_samples}, preserve convexity of the objective. In the convex setting, the projected
stochastic subgradient method \eqref{eq:algo_sgd_beta} admits the following bound on the averaged objective values.

\begin{proposition}
\label{prop:rates_cvx}
Under \Cref{ass:bounds_noise_beta}, let $(w_k)_{k\in\mathbb{N}}$ be generated by
\eqref{eq:algo_sgd_beta}. Assume that all functions $\widetilde H_i^k$, $H_i$ for 
$i=1,\ldots,n$, $k\in\mathbb{N}$, are convex. Let $F^*:=\min_K F$ and $w^*\in\argmin_K F$. Then, for every $N\geq 1$,
\[
\begin{aligned}
    \mathbb{E}\bigl[F(w_{R_N})-F^*\bigr]
    & \leq
    \frac{\|w_0 - w^*\|^2}{2A_N}
+
    B_N^{\mathrm{cvx}}
    \\ &+
    \frac{2C_K}{A_N}
    \sum_{k=0}^{N-1}
    \alpha_k\beta_k|\log\beta_k|
    +
    \frac{\widetilde M_K^2}{2A_N}
    \sum_{k=0}^{N-1}\alpha_k^2,
\end{aligned}
\]
where $B_N^{\mathrm{cvx}}
    :=
    \frac{1}{A_N}
    \sum_{k=0}^{N-1}
    \alpha_k
    \mathbb{E}
    \left[
        \bar b_k(w_k)+\bar b_k(w^*)
    \right]$, $A_N:=\sum_{k=0}^{N-1}\alpha_k$ and
$\mathbb{P}(R_N=k)=\alpha_k/A_N$ for $k=0,\ldots,N-1$. In particular, if
$\alpha_k=\alpha_0/\sqrt{k+1}$ and $\beta_k=\beta_0/(k+1)$, with
$\alpha_0,\beta_0>0$, then
\[
    \mathbb{E}\bigl[F(w_{R_N})-F^*\bigr]
    =
    O\left(
        \frac{\log N}{\sqrt N}
        +
        B_N^{\mathrm{cvx}}
    \right),
\]
\end{proposition}

\begin{proof}
By nonexpansiveness of the projection,
	\[
	\begin{aligned}
		\|w_{k+1}-w^*\|^2
		\leq{}&
		\|w_k-w^*\|^2
		-2\alpha_k
		\left\langle
		\nabla\widetilde H_{i_k}^k(w_k),w_k-w^*
		\right\rangle \\
		&+
		\alpha_k^2
		\left\|
		\nabla\widetilde H_{i_k}^k(w_k)
		\right\|^2.
	\end{aligned}
	\]
	Since $\widetilde H_{i_k}^k$ is convex, $\left\langle
	\nabla\widetilde H_{i_k}^k(w_k),w_k-w^*
	\right\rangle
	\geq
	\widetilde H_{i_k}^k(w_k)-\widetilde H_{i_k}^k(w^*).$
	Taking conditional expectations and using the uniform sampling of $i_k$
	together with the approximation bound from \Cref{ass:bounds_noise_beta} and \eqref{eq:b_bar} yields
	\[
	\begin{aligned}
		\mathbb E\left[
		\widetilde H_{i_k}^k(w_k)-\widetilde H_{i_k}^k(w^*)
		\,\middle|\,w_k
		\right]
		\geq{}&
		F^{\beta_k}(w_k)-F^{\beta_k}(w^*) \\
		&-\bar b_k(w_k)-\bar b_k(w^*).
	\end{aligned}
	\]
	Moreover, the regularization-error bound gives
	\[
	F^{\beta_k}(w_k)-F^{\beta_k}(w^*)
	\geq
	F(w_k)-F(w^*)
	-2C_K\beta_k|\log\beta_k|.
	\]
	Therefore,
	\[
	\begin{aligned}
		2\alpha_k\bigl(F(w_k)-F^*\bigr)
		\leq{}&
		\|w_k-w^*\|^2
		-
		\mathbb E\left[
		\|w_{k+1}-w^*\|^2
		\,\middle|\,w_k
		\right] \\
		&+
		2\alpha_k
		\bigl(
		\bar b_k(w_k)+\bar b_k(w^*)
		\bigr) \\
		&+
		4C_K\alpha_k\beta_k|\log\beta_k|
		+
		\alpha_k^2\widetilde M_K^2.
	\end{aligned}
	\]
	Taking expectations and summing from $k=0$ to $N-1$ gives
	\[
	\begin{aligned}
		\frac{1}{A_N}
		\sum_{k=0}^{N-1}
		\alpha_k\,
		\mathbb E\bigl[F(w_k)-F^*\bigr]
		\leq{}&
		\frac{\|w_0-w^*\|^2}{2A_N}
		+B_N^{\mathrm{cvx}} \\
		&+
		\frac{2C_K}{A_N}
		\sum_{k=0}^{N-1}
		\alpha_k\beta_k|\log\beta_k|
		+
		\frac{\widetilde M_K^2}{2A_N}
		\sum_{k=0}^{N-1}\alpha_k^2.
	\end{aligned}
	\]
	Using the
	definition of $R_N$ proves the general bound. Finally, the stated
	rate follows from the standard estimates for the corresponding sums.
\end{proof}

\subsection{Weakly convex case}

We now consider the weakly convex setting. The stationarity measure is the
gradient of the \emph{Moreau envelope} of the objective, restricted to the compact convex set $K \subset \R^p$. For
$\tau>0$, define
\[
    \varphi_{1/\tau}(w)
    :=
    \min_{u\in K}
    \left\{
        F(u)+\frac{\tau}{2}\|u-w\|^2
    \right\},\]
and $\operatorname{prox}_{\tau^{-1}F}(w)
    :=
    \argmin_{u\in K}
    \left\{
        F(u)+\frac{\tau}{2}\|u-w\|^2
    \right\}.$
If $F$ is $\tau_0$-weakly convex, then for every $\tau>\tau_0$,
$\varphi_{1/\tau}$ is continuously differentiable \cite{moreau1965proximite}.

\begin{proposition}
\label{prop:rates_weakly_cvx}
Under \Cref{ass:bounds_noise_beta}, let $(w_k)_{k\in\mathbb{N}}$ be generated by
\eqref{eq:algo_sgd_beta}. Assume there exists $\tau_0>0$ such that all functions $\widetilde H_i^k$, $H_i$ for 
$i=1,\ldots,n$ and $k\in\mathbb{N}$, are $\tau_0$-weakly convex. Let $\tau>\tau_0$, and for each $k \geq 0$, let $\widehat w_k:=\operatorname{prox}_{\tau^{-1}F}(w_k),$
and let $\varphi^*:=\inf_{w\in K}F(w)$. Then, for every $N\geq 1$,
\[
\begin{aligned}
    \mathbb{E}
    \left[
        \|\nabla\varphi_{1/\tau}(w_{R_N})\|^2
    \right]
    \leq
    \frac{2\tau(\varphi_{1/\tau}(w_0)-\varphi^*)}
    {(\tau-\tau_0)A_N}
    &+
    \frac{2\tau^2}{(\tau-\tau_0)} B_N^{\mathrm{wc}}
    \\
    +
    \frac{4C_K\tau^2}{(\tau-\tau_0)A_N}
    \sum_{k=0}^{N-1}
    \alpha_k\beta_k|\log\beta_k|
    &+
    \frac{\tau\widetilde M_K^2}
    {(\tau-\tau_0)A_N}
    \sum_{k=0}^{N-1}\alpha_k^2,
\end{aligned}
\]
where given $\bar{b}_k$ defined by \eqref{eq:b_bar}, $B_N^{\mathrm{wc}}
    :=
    \frac{1}{A_N}
    \sum_{k=0}^{N-1}
    \alpha_k
    \mathbb{E}
    \left[
        \bar b_k(w_k)+\bar b_k(\widehat w_k)
    \right]$,    $A_N:=\sum_{k=0}^{N-1}\alpha_k$ and  $\mathbb{P}(R_N=k)=\alpha_k/A_N$ for $k=0,\ldots,N-1$. In particular, if
$\alpha_k=\alpha_0/\sqrt{k+1}$ and $\beta_k=\beta_0/(k+1)$, with
$\alpha_0,\beta_0>0$, then
\[
    \mathbb{E}
    \left[
        \|\nabla\varphi_{1/\tau}(w_{R_N})\|^2
    \right]
    =
    O\left(
        \frac{\log N}{\sqrt N}
        +
        B_N^{\mathrm{wc}}
    \right),
\]

\end{proposition}

\begin{proof}
	By nonexpansiveness of the projection,
	\[
	\begin{aligned}
		\|w_{k+1}-\widehat w_k\|^2
		\leq{}&
		\|w_k-\widehat w_k\|^2
		-2\alpha_k
		\left\langle
		\nabla\widetilde H_{i_k}^k(w_k),
		w_k-\widehat w_k
		\right\rangle \\
		&+
		\alpha_k^2
		\left\|
		\nabla\widetilde H_{i_k}^k(w_k)
		\right\|^2.
	\end{aligned}
	\]
	Since $\widetilde H_{i_k}^k$ is $\tau_0$-weakly convex,
	\[
	\left\langle
	\nabla\widetilde H_{i_k}^k(w_k),
	w_k-\widehat w_k
	\right\rangle
	\geq
	\widetilde H_{i_k}^k(w_k)
	-\widetilde H_{i_k}^k(\widehat w_k)
	-\frac{\tau_0}{2}
	\|w_k-\widehat w_k\|^2.
	\]
	Taking conditional expectations and using the uniform sampling of $i_k$
	together with \Cref{ass:bounds_noise_beta} and the averaged error \eqref{eq:b_bar}, we obtain
	\[
	\begin{aligned}
		\mathbb E\left[
		\widetilde H_{i_k}^k(w_k)
		-\widetilde H_{i_k}^k(\widehat w_k)
		\,\middle|\,w_k
		\right]
		\geq{}&
		F^{\beta_k}(w_k)-F^{\beta_k}(\widehat w_k) \\
		&-\bar b_k(w_k)-\bar b_k(\widehat w_k).
	\end{aligned}
	\]
	Moreover, $F^{\beta_k}(w_k)-F^{\beta_k}(\widehat w_k)
	\geq
	F(w_k)-F(\widehat w_k)
	-2C_K\beta_k|\log\beta_k|$ and by the definition of $\widehat w_k$, $F(w_k)-F(\widehat w_k)
	\geq
	\frac{\tau}{2}\|w_k-\widehat w_k\|^2.$
	Consequently,
	\[
	\begin{aligned}
		\mathbb E\left[
		\|w_{k+1}-\widehat w_k\|^2
		\,\middle|\,w_k
		\right]
		\leq{}&
		\|w_k-\widehat w_k\|^2
		-\alpha_k(\tau-\tau_0)
		\|w_k-\widehat w_k\|^2 \\
		&+
		2\alpha_k\left(
		\bar b_k(w_k)+\bar b_k(\widehat w_k)
		+2C_K\beta_k|\log\beta_k|
		\right)
		+\alpha_k^2\widetilde M_K^2.
	\end{aligned}
	\]
	Using $\varphi_{1/\tau}(w_k)
	=
	F(\widehat w_k)
	+\frac{\tau}{2}\|w_k-\widehat w_k\|^2$
	and $\varphi_{1/\tau}(w_{k+1})
	\leq
	F(\widehat w_k)
	+\frac{\tau}{2}\|w_{k+1}-\widehat w_k\|^2,$
	we deduce
	\[
	\begin{aligned}
		&\frac{\alpha_k\tau(\tau-\tau_0)}{2}
		\|w_k-\widehat w_k\|^2
		\leq{}
		\varphi_{1/\tau}(w_k)
		-
		\mathbb E\left[
		\varphi_{1/\tau}(w_{k+1})
		\,\middle|\,w_k
		\right] \\
		&+
		\alpha_k\tau\left(
		\bar b_k(w_k)+\bar b_k(\widehat w_k)
		+2C_K\beta_k|\log\beta_k|
		\right)
		+\frac{\alpha_k^2\tau}{2}\widetilde M_K^2.
	\end{aligned}
	\]
	Since $\nabla\varphi_{1/\tau}(w_k)
	=
	\tau(w_k-\widehat w_k),$
	it follows that
	\[
	\begin{aligned}
		& \frac{\alpha_k}{2}
		\left(1-\frac{\tau_0}{\tau}\right)
		\|\nabla\varphi_{1/\tau}(w_k)\|^2
		\leq{}
		\varphi_{1/\tau}(w_k)
		-
		\mathbb E\left[
		\varphi_{1/\tau}(w_{k+1})
		\,\middle|\,w_k
		\right] \\
		&+
		\alpha_k\tau\left(
		\bar b_k(w_k)+\bar b_k(\widehat w_k)
		+2C_K\beta_k|\log\beta_k|
		\right)
		+\frac{\alpha_k^2\tau}{2}\widetilde M_K^2.
	\end{aligned}
	\]
	Taking expectations, summing from $k=0$ to $N-1$, and using
	$\varphi_{1/\tau}(w_N)\geq\varphi^*$ gives
	\[
	\begin{aligned}
		& \frac{1}{A_N}
		\sum_{k=0}^{N-1}
		\alpha_k\,
		\mathbb E\left[
		\|\nabla\varphi_{1/\tau}(w_k)\|^2
		\right]
		\leq{}
		\frac{2\tau\bigl(
			\varphi_{1/\tau}(w_0)-\varphi^*
			\bigr)}
		{(\tau-\tau_0)A_N} 
		+
		\frac{2\tau^2}{\tau-\tau_0}
		B_N^{\mathrm{wc}} \\
		& +
		\frac{4C_K\tau^2}{(\tau-\tau_0)A_N}
		\sum_{k=0}^{N-1}
		\alpha_k\beta_k|\log\beta_k| 
		+
		\frac{\tau^2\widetilde M_K^2}
		{(\tau-\tau_0)A_N}
		\sum_{k=0}^{N-1}\alpha_k^2.
	\end{aligned}
	\]
	The definition of $R_N$ yields the claimed general estimate. The stated rate
	then follows from the standard estimates for the corresponding sums.
\end{proof}
\section{Application to WDRO}
\label{sec:proofs}
We now prove the main results of \Cref{sec:main_results}. They combine approximation results from \Cref{section:convergence_entropic_general} and convergence rates from \Cref{sec:convergence_rates_general}.

\subsection{Gradient  and critical set convergence}
\label{subsection:convergence_gradient_proof}

\begin{proof}[of \Cref{th:gradient_conv_ell_c}] We prove the qualitative convergence result. For each $1 \leq i \leq n$, we set   $h_i(\theta, \lambda,z) =   \ell(\theta,z) + \lambda (\rho - c(\xi_i, z))$ and  $H_i(\theta,\lambda) = \max_{z \in \mathcal{Z}} h_i(\theta,\lambda,z)$, which is well defined for any $\lambda \geq \lambda_{\min}$ under  \Cref{ass:regularity_ell_c_merged}.3 and we denote the corresponding entropic regularization by $H_i^{\beta}$. The abstract qualitative results in
	\Cref{th:graphical_convergence_one_batch,cor:mini_batch_graphical_convergence}
	apply locally on an open neighborhood of $K$, since their proofs use
	only local bounds near $K$. Under \Cref{ass:regularity_ell_c_merged},
	with the regularity and growth conditions imposed on a neighborhood
	of $\Theta$, the strict inequality $\bar\lambda<\lambda_{\min}$
	ensures that \Cref{ass:regularity_h} and \Cref{ass:differentiable_h} hold locally
	near $K=\Theta\times\Lambda$ for each $h_i$. Applying  \Cref{cor:mini_batch_graphical_convergence} then gives the first claim. The second claim, concerning convergence of the critical sets, follows from \Cref{lem:graphconv_to_critconv}.
\end{proof}

\begin{proof}[of \Cref{lem:LLN_H_beta_m}]
	We fix $\beta>0$ and $i\in\{1,\ldots,n\}$. We set $w=(\theta,\lambda)$ and
	\begin{equation*}
		\begin{aligned}
			u_i(w,z) := \exp\!\left(
			\frac{\ell(\theta,z)+\lambda(\rho-c(\xi_i,z))}{\beta}
			\right), \quad q_i(w,z) := \begin{pmatrix}
				\nabla_\theta\ell(\theta,z)\\
				\rho-c(\xi_i,z)
			\end{pmatrix}.
		\end{aligned}
	\end{equation*}
	We define the quantities
	\begin{align*}
		& A(w):=\E_{z \sim \pi_0}[u_i(w,z)],
		& A_m(w):=\frac1m\sum_{j=1}^m u_i(w,z_j),\\
		& B(w):=\E_{z \sim \pi_0}[u_i(w,z)q_i(w,z)], 
		& B_m(w):=\frac1m\sum_{j=1}^m u_i(w,z_j)q_i(w,z_j).
	\end{align*}
	With these notations,
	\begin{equation}\label{eq:gradient_H_ratio}
		\begin{aligned}
			H_i^\beta &= \beta\log A,
			&\qquad H_i^{\beta,m} &= \beta\log A_m,\\
			\nabla H_i^\beta &= \frac{B}{A},
			&\qquad \nabla H_i^{\beta,m} &= \frac{B_m}{A_m}.
		\end{aligned}
	\end{equation}
	We first verify the integrability needed to apply the uniform law of large numbers \cite[Th. 9.60]{shapiro2021lectures}.
	Let $a:=\lambda_{\min}-\bar\lambda>0$. By
	\Cref{ass:regularity_ell_c_merged}, for every $w\in\Theta\times\Lambda$, $\ell(\theta,z)+\lambda(\rho-c(\xi_i,z))
	\leq L+\lambda_{\max}\rho-a c(\xi_i,z).$
	Consequently, with $C_\beta:=\exp((L+\lambda_{\max}\rho)/\beta)$,
	\begin{equation*}
		0<u_i(w,z)\leq C_\beta \exp(-a c(\xi_i,z)/\beta)\leq C_\beta.
	\end{equation*}
	Moreover,
	\begin{align*}
		\|u_i(w,z)q_i(w,z)\|
		&\leq C_\beta \exp(-a c(\xi_i,z)/\beta)
		\bigl(\varphi(z)+\rho+c(\xi_i,z)\bigr)\\
		&\leq C_\beta\left(\varphi(z)+\rho+\frac{\beta}{a\exp(1)}\right),
	\end{align*}
	where we used $t \exp(-at/\beta)\leq\beta/(a\exp(1))$ for $t\geq0$.
	The right-hand side is $\pi_0$-integrable. Since the two integrands are
	continuous in $w$ and $\Theta\times\Lambda$ is compact, the uniform law of
	large numbers yields, almost surely,
	\begin{equation}\label{eq:uniform_AB}
		\|A_m-A\|_\infty\longrightarrow0
		\quad \text{and} \quad 
		\|B_m-B\|_\infty\longrightarrow0,
	\end{equation}
	where the suprema $\|\cdot\|_\infty$ are taken over $\Theta\times\Lambda$. It remains to pass these limits through the logarithm and the ratio in
	\eqref{eq:gradient_H_ratio}. The dominated convergence theorem shows that $A$ and $B$
	are continuous on $\Theta\times\Lambda$. In addition, since $A$ is positive and continuous on the compact set $\Theta\times\Lambda$, $\delta_i:=\inf_{w\in\Theta\times\Lambda}A(w)>0.$
	
	By \eqref{eq:uniform_AB}, almost surely $\inf_w A_m(w)\geq\delta_i/2$ for
	all sufficiently large $m$. 
	Elementary estimates then give the desired convergence, that is
	\begin{equation*}
		\|\beta\log A_m-\beta\log A\|_\infty
		\leq \frac{2\beta}{\delta_i}\|A_m-A\|_\infty\longrightarrow0
	\end{equation*}
	and $\left\|\frac{B_m}{A_m}-\frac{B}{A}\right\|_\infty
	\leq \frac{2}{\delta_i}\|B_m-B\|_\infty
	+\frac{2\|B\|_\infty}{\delta_i^2}\|A_m-A\|_\infty
	\longrightarrow0,$
	where the suprema are taken over $\Theta \times \Lambda$. Both suprema converge to zero almost surely. Since there are
	only finitely many indices $i$, their probability-one events may be
	intersected and their eventual sample-size thresholds maximized to obtain the desired $\varepsilon$ threshold and the minimal sample size $m_{\beta,\varepsilon}$.

For the expectation claim, set $Z \sim \pi_0$, $K=\Theta\times\Lambda$,
$U=L+\lambda_{\max}\rho$ and
$B(z)=a(z)+\lambda_{\max}(\rho+\max_i c(\xi_i,z))$.
Then $B$ is integrable and $-B(z)\leq h_i(w,z)\leq U$ uniformly in $i,w$.
Jensen's inequality gives
\[
-\frac1m\sum_{j=1}^m B(z_j)\leq H_i^{\beta,m}(w)\leq U,
\qquad -\mathbb E B(Z)\leq H_i^\beta(w)\leq U.
\]
Hence the uniform objective error satisfies
\[
\max_{1\leq i\leq n}\sup_{w\in K}
\left|H_i^{\beta,m}(w)-H_i^\beta(w)\right|
\leq
2U+\frac1m\sum_{j=1}^m B(z_j)+\mathbb E B(Z).
\]
The empirical averages of the integrable random variable $B(Z)$ converge
in $L^1$ to $\E [B(Z)]$ and are therefore uniformly integrable. The almost-sure uniform convergence proved above thus also holds in $L^1$ by Vitali's convergence theorem, yielding the claimed
deterministic threshold $m_{\beta,\varepsilon}$.
\end{proof}

\begin{proof}[of \Cref{th:quantitative_approx_graph}] We prove the quantitative convergence result. Our goal is to apply \Cref{prop:quantitative_rate_h} with  the functions $h_i(\theta,\lambda,z) = \ell(\theta,z) + \lambda (\rho - \|\xi_i - z\|^2)$ for $i =1, \ldots, n$ and the compact set $K = \Theta \times \Lambda$.

	We fix $i\in\{1,\dots,n\}$ and we verify that \Cref{ass:L_smooth_h} and \Cref{ass:error_bound} hold on $K$ for $h := h_i$, $\xi := \xi_i$, under \Cref{ass:WDRO_smooth_quadratic}. \Cref{ass:L_smooth_h}.1 corresponds to  \Cref{ass:WDRO_smooth_quadratic}.1 and \Cref{ass:L_smooth_h}.2 is implied by \Cref{ass:WDRO_smooth_quadratic}.4.

In order to verify \Cref{ass:error_bound}, we first show that the
maximizers of $h(\theta,\lambda,\cdot)$, denoted by $M(\theta,\lambda)$, are contained
in a common compact set, uniformly over
$(\theta,\lambda)\in K$. Under
\Cref{ass:WDRO_smooth_quadratic}.4, we have
$\ell(\theta,z)\leq L+\bar\lambda\|z-\xi\|^2$ where
$\bar\lambda\in[0,\lambda_{\min})$. Hence, for every
$\lambda\in\Lambda$,
\begin{equation}
	\label{eq:right-hand-side}
	\begin{aligned}
		h(\theta,\lambda,z)
		&=
		\ell(\theta,z)+\lambda\bigl(\rho-\|z-\xi\|^2\bigr) \\
		&\leq
		L+\lambda_{\max}\rho
		-(\lambda_{\min}-\bar\lambda)\|z-\xi\|^2.
	\end{aligned}
\end{equation}
Since $\lambda_{\min}-\bar\lambda>0$, the right-hand side is
anti-coercive. Together with continuity, this implies that
$H(\theta,\lambda):=\max_{z\in\mathcal Z}
h(\theta,\lambda,z)$ is well defined and continuous. Thus,
$\underline H_{K}:=\min_{i=1,\ldots,n}\min_{(\theta,\lambda)\in K}
H_i(\theta,\lambda)$ is finite and for every
$m\in M(\theta,\lambda)$, $\underline H_{K}
\leq
L+\lambda_{\max}\rho
-(\lambda_{\min}-\bar\lambda)\|m-\xi\|^2.$
Therefore, setting
$R_{K}^2:=(L+\lambda_{\max}\rho-\underline H_{K})/
(\lambda_{\min}-\bar\lambda)$, all the maximizers are contained in
the compact set
$\Xi_0:=\bigcup_{i=1}^n\closedball(\xi_i,R_{K})$. Let $\Xi$ be the closed neighborhood of $\Xi_0$ with radius $1$, which is also compact. This enlargement prevents the maximizers from lying on the boundary of the compact set.

We now verify \Cref{ass:error_bound}.1. By the joint Lipschitz
continuity of $\nabla_z\ell$ on the compact set
$\Theta\times\Xi$, as stated in
\Cref{ass:WDRO_smooth_quadratic}.3, there exists $L_\Xi>0$ such that
the descent lemma holds: for every $m,z\in\Xi$,
\[
h(\theta,\lambda,z)
\geq
h(\theta,\lambda,m)
+
\left\langle
\nabla_z h(\theta,\lambda,m),z-m
\right\rangle
-
\frac{L_\Xi+2\lambda_{\max}}{2}\|z-m\|^2.
\]
For $m\in M(\theta,\lambda)$, we have
$\nabla_z h(\theta,\lambda,m)=0$ and
$h(\theta,\lambda,m)=H(\theta,\lambda)$. Choosing
$m\in M(\theta,\lambda)$ such that
$\|z-m\|=\dist(z,M(\theta,\lambda))$, we obtain
\[
h(\theta,\lambda,z)
\geq
H(\theta,\lambda)
-
\frac{L_\Xi+2\lambda_{\max}}{2}
\dist\bigl(z,M(\theta,\lambda)\bigr)^2.
\]
Hence, \Cref{ass:error_bound}.1 holds with $q=2$ and
$c_K=(L_\Xi+2\lambda_{\max})/2$ whenever
$\dist(z,M(\theta,\lambda))\leq1$.

We next verify \Cref{ass:error_bound}.2. Applying
\Cref{ass:WDRO_smooth_quadratic}.5 to the compact set $\Xi$, we
obtain $a_K>0$ and $\bar r_K>0$ such that, for every
$i=1,\ldots,n$, $(\theta,\lambda)\in K$,
$m\in M(\theta,\lambda)$, and $r\in(0,\bar r_K]$,
$\pi_0(\closedball(m,r))\geq a_Kr^d$. Therefore, taking
$r_K:=\min\{1,\bar r_K\}$ gives both parts of
\Cref{ass:error_bound} with the same radius.

Finally, let $\tau_\ell>0$ be given by
\Cref{ass:WDRO_smooth_quadratic}.2, so that $\ell(\cdot,z)$ is
$\tau_\ell$-weakly convex for every $z\in\mathcal Z$. Since
$\lambda\mapsto\lambda(\rho-\|z-\xi\|^2)$ is affine,
$(\theta,\lambda)\mapsto h(\theta,\lambda,z)$ is
$\tau_\ell$-weakly convex. Thus, the desired weak convexity property holds
with $\tau_h=\tau_\ell$, uniformly in $z$. Applying \Cref{prop:quantitative_rate_h}on $K$ gives the desired estimates.

\end{proof}
\subsection{Stochastic gradient methods}
\label{subsection:SGD_proofs}
We now prove \Cref{prop:convergence_SGD_vanish_reg} for vanishing steps and regularization parameters.

\begin{proof}[of \Cref{prop:convergence_SGD_vanish_reg}] This is an application of the more general results, \Cref{prop:rates_cvx} and \Cref{prop:rates_weakly_cvx} with the choice $\widetilde{H}_i^k = H^{\beta_k,m_k}_i$, where for any $\beta >0$ and sample size $m$, $H_i^{\beta,m}$ is defined as \eqref{eq:H_beta_m}.

	Let $\Delta_m$ denote the probability simplex on $\{1,\ldots,m\}$, and let
	$q^{(m)}=(1/m,\ldots,1/m)$ be the uniform distribution on these indices.
	The Donsker--Varadhan variational formula \cite{donsker1983asymptotic} gives
	\[
	H_i^{\beta,m}(\theta,\lambda)
	=
	\sup_{u\in\Delta_m}
	\left\{
	\sum_{j=1}^m u_j
	\bigl[\ell(\theta,z_j)+\lambda\bigl(\rho-c(\xi_i,z_j)\bigr)\bigr]
	- \beta\,\KL(u\|q^{(m)})
	\right\}.
	\]
	Here $\KL(u\|q^{(m)})=\sum_{j=1}^m u_j\log(m u_j)$ is the
	Kullback--Leibler divergence, with the convention $0\log 0=0$.
	If $\ell(\cdot,z)$ is $\tau_\ell$-weakly convex uniformly in $z$, then,
	for every $u\in\Delta_m$, the function inside the supremum is
	$\tau_\ell$-weakly convex in $w$: the transport cost term is affine in
	$\lambda$, while the coefficients $u_j$ are nonnegative and sum to
	one. Since the pointwise supremum of functions sharing the same weak
	convexity constant is again $\tau_\ell$-weakly convex,
	$H_i^{\beta,m}$ is $\tau_\ell$-weakly convex. The same argument with
	$\tau_\ell=0$ shows that $H_i^{\beta,m}$ is convex whenever
	$\ell(\cdot,z)$ is convex.

We now verify the three conditions in \Cref{ass:bounds_noise_beta}. Condition 1 follows from \eqref{eq:m_k_bound} with $b_{i,k}(w)=\bar b$. The same bound holds conditionally on the history, including at the current iterate and its proximal point. Condition 3 is a consequence of the quantitative approximation result, \Cref{th:quantitative_approx_graph} which holds under \Cref{ass:WDRO_smooth_quadratic}. It remains to verify Condition 2 under \Cref{ass:WDRO_smooth_quadratic}. We use the shorthand $w=(\theta,\lambda)$ throughout. For every $i=1,\ldots,n$, let
\[
h_i(w,z)
:=
\ell(\theta,z)
+\lambda\bigl(\rho-c(\xi_i,z)\bigr).
\]
The gradient of the sampling approximation can be written with the empirical Gibbs measure
\[
\nabla H_i^{\beta,m}(w)
=
\int_{\mathcal Z}
\nabla_w h_i(w,z)\,
\di\pi_{\beta,w}^{m}(z), \qquad \pi_{\beta,w}^{m} \propto \sum_{j=1}^m e^{h_i(w,z_j)/\beta} \delta_{z_j}
\]
where $z_1, \ldots, z_m$ are i.i.d.\ from $\pi_0$. By Jensen's inequality,
\begin{equation}
	\label{eq:sample_gradient_Jensen}
	\left\|\nabla H_i^{\beta,m}(w)\right\|^2
	\leq
	\int_{\mathcal Z}
	\left\|\nabla_w h_i(w,z)\right\|^2
	\,\di\pi_{\beta,w}^{m}(z).
\end{equation}
We easily verify that $\nabla_w h_i(w,z)
=
\begin{pmatrix}
	\nabla_\theta\ell(\theta,z)\\
	\rho-c(\xi_i,z)
\end{pmatrix}$ and 
$$\|\nabla_w h_i(w,z)\|^2 = \|\nabla_\theta \ell (\theta, z)\|^2 + (\rho - \|\xi_i - z\|^2)^2$$
By item 4, there exists $C_K>0$ such that $\left\|\nabla_\theta \ell(\theta,z)\right\|^2
+ \|z\|^4\leq C_K\bigl(1+\Psi(z)\bigr)$ for every $w\in K$ and $z\in\mathcal Z$. Thus there exist $b >0$ and $a > 0$ such that for all $z \in \mathcal{Z}$, $\|\nabla_w h_i(w,z)\|^2 \leq b(a + \Psi(z)).$ Using this inequality in \eqref{eq:sample_gradient_Jensen} gives
\begin{equation}
	\label{eq:sample_gradient_Psi_bound}
	\left\|\nabla H_i^{\beta,m}(w)\right\|^2
	\leq
	b \left(
	a + 
	\int_{\mathcal Z}\Psi(z)\,
	\di\pi_{\beta,w}^{m}(z)
	\right).
\end{equation}

This holds for any $w \in K$ and every $\beta >0$. Now, under items 4 and 5 of \Cref{ass:WDRO_smooth_quadratic}, we can apply  \Cref{lem:sufficient_uniform_moment_empirical_gibbs} to obtain the desired bound
\begin{equation*}
	\sup_{i = 1, \ldots,n} \sup_{w \in K, \beta > 0} \E \left[\|\nabla H_i^{\beta,m}(w)\|^2 \right] := M_K< \infty
\end{equation*}
Since the preceding estimate is uniform in
$i$, $w\in K$, $\beta>0$, and $m\geq1$, and the batch at iteration $k$ is conditionally i.i.d.\ with law $\pi_0$ given $w_k$, we
obtain
\[
\sup_{i=1,\ldots,n} \sup_{k\in\mathbb N}
\mathbb E\left[
\left\|\nabla H_i^{\beta_k,m_k}(w_k)\right\|^2
\,\middle|\,w_k
\right]
\leq M_K
\quad\text{a.s.}
\]
\end{proof} 

\section{Conclusion}
We studied approximation guarantees for entropic smoothing in Wasserstein distributionally robust optimization. We established qualitative and quantitative convergence results for sampled gradient oracles and their implications for stochastic gradient methods. We considered both fixed regularization, with approximation guarantees on compact parameter sets, and varying regularization, with guarantees along algorithmic trajectories.

Future work includes developing sampling methods for the Gibbs distribution and evaluating the proposed approaches on learning tasks. Vanishing regularization may pose numerical challenges, motivating adaptive schedules and theoretical guarantees on the smoothing levels needed for good learning performance.

\bibliographystyle{plain}
\bibliography{references}

\end{document}